\documentclass[a4paper,preprint,12pt]{elsarticle}
\usepackage[english]{babel}
\usepackage[latin1]{inputenc}
\usepackage{amssymb}
\usepackage{amsmath}         
\usepackage{amsthm}             
\usepackage{pb-diagram}              
\usepackage{epstopdf}      
\usepackage{graphicx}             
 \usepackage{epsfig}       
\usepackage{color}       
\usepackage{fancyhdr}   
 \usepackage{natbib}      
\usepackage{rotating}
\usepackage[T1]{fontenc}        
\usepackage{url}                
\usepackage{mathtools}                 
  
\usepackage{vmargin}            
\usepackage[title]{appendix} 

\usepackage{lineno}
\usepackage{setspace}                      
\biboptions{square} 
\usepackage[colorlinks=true,linkcolor= blue,urlcolor=blue,citecolor=blue]{hyperref} 
\usepackage[ruled,vlined]{algorithm2e} 

\newcommand{\hh}{h}

 \newcommand{\bo}[1]{\mathbf{#1}}  
\def\indic{\hbox{1\kern-.24em\hbox{I}}}      

\newcommand{\var}{\mathbb{V}}  
\newcommand{\esp}{\mathbb{E}}

\newcommand{\M}{f}

\newcommand{\T}{T}    
\newcommand{\N}{\mathbb{N}}      
\newcommand{\NN}{n}

\newcommand{\norme}[1]{\left|\left| #1 \right|\right|_{2}}

\newcommand{\norminf}[1]{\left|\left| #1 \right|\right|_{\infty}}

\newcommand{\rr}{r}

\newcommand{\x}{x}
\newcommand{\X}{X}

\newcommand{\R}{\mathbb{R}}

\newtheorem{prop}{Proposition}{\bf}{\it} 
\newtheorem{defi}{Definition}{\bf}{\it}  

\newtheorem{theorem}{Theorem}{\bf}{\it}     
\newtheorem{assump}{Assumption}{\bf}{\it}  
\newtheorem{lemma}{Lemma}{\bf}{\it}          
\newtheorem{rem}{Remark}{\bf}{\it} 
\newtheorem{corollary}{Corollary}{\bf}{\it} 

\catcode`\"=\active
\catcode`\"=\active
\def"{\og\ignorespaces}
\def"{{\fg}}

\makeatletter
\def\ps@pprintTitle{%
  \let\@oddhead\@empty
  \let\@evenhead\@empty
  \def\@oddfoot{\reset@font\hfil\thepage\hfil}
  \let\@evenfoot\@oddfoot
}
\makeatother

\begin{document}  
\begin{frontmatter}           
\title{Optimal ANOVA-based emulators of models with(out) derivatives}
\author[a,b]{Matieyendou Lamboni\footnote{Corresponding author: matieyendou.lamboni[at]gmail.com or [at]univ-guyane.fr; Freb 26, 2025}}              
\address[a]{University of Guyane, Department DFR-ST, 97346 Cayenne, French Guiana, France}
\address[b]{228-UMR Espace-Dev, University of Guyane, University of R\'eunion, IRD, University of Montpellier, France.}                                      
                                                           
\begin{abstract}     
This paper proposes new ANOVA-based approximations of functions and emulators of high-dimensional models using either available derivatives or local stochastic evaluations of such models. Our approach makes use of sensitivity indices to design adequate structures of emulators. For high-dimensional models with available derivatives, our derivative-based emulators reach dimension-free mean squared errors (MSEs) and parametric rate of convergence (i.e., $\mathsf{O}(N^{-1})$). This approach is extended to cope with every  model (without available derivatives) by deriving global emulators that account for the local properties of models or simulators. Such generic emulators enjoy dimension-free biases, parametric rates of convergence and MSEs that depend on the dimensionality. Dimension-free MSEs are obtained for high-dimensional models with particular inputs' distributions. Our emulators are also competitive in dealing with different distributions of the input variables and for selecting inputs and interactions. Simulations show the efficiency of our approach.
      
\end{abstract}     
\end{frontmatter}                

       
\section{Introduction}    
 Derivatives are sometime available in modeling either by the nature of observations of the phenomena of interest (\cite{morris93,solak02} and the reference therein) or by low-cost evaluations of exact derivatives for some classes of PDE/ODE-based models thanks to adjoint methods (\cite{dimet86,dimet97,cacuci05,gunzburger03,borzi12,ghanem17}).   
Models defined via their rates of change w.r.t. their inputs; implicit functions (defined via their derivatives) are instances. Additionally and particularly for complex models or simulators, efficient estimators of gradients, second-order derivatives using stochastic approximations have been provided in \cite{agarwal10,bach16,akhavan20,lamboni24axioms,lamboni24stats}. Being able to reconstruct functions using such derivatives is worth investigating (\cite{chkifa13}), as having a practical, fast-evaluated model tnat links stochastic parameters and/or stochastic initial conditions to the output of interest of PDE/ODE models remains a challenge due to numerous uncertain inputs (\cite{patil23}). \\   
Moreover, the first-order derivatives of models are used for quickly selecting non-relevant input variables of simulators or functions, leading to effective screening measures. Efficient variance-based screening methods rely on the upper-bounds of generalized sensitivity indices, including Sobol' indices (see \cite{sobol09,kucherenko09,lamboni13,roustant14,roustant17,lamboni20,lamboni22} for independent inputs and \cite{lamboni21} for non-independent variables). \\       
                
For high-dimensional simulators, dimension reductions via screening input variables are often performed before building their emulators using Gaussian processes (kriging models) (\cite{krige51,currin91,oakley04,conti10,kennedy01}), polynomial chaos expansions (\cite{dongbin02,ghanem91,sudret08}), SS-ANOVA (\cite{wahba90,wong16}) or other machine learning approaches (\cite{friedman08,horiguchi24}). Indeed, such emulators rely on nonparametric or semi-parametric regressions and struggle to reconstruct simulators for a moderate to large number of inputs (e.g., \cite{horiguchi24}). Often, non-parametric rates of convergence are achieved by such emulators (see \cite{wong16} for SS-ANOVA and \cite{migliorati14,hampton15,cohen18} for polynomial chaos expansions). Regarding the stability and accuracy of polynomial chaos expansions, the number of model runs needed is firstly estimated at  the square of the dimension of the basis used  (\cite{migliorati14}), and then reduce at that dimension up to a logarithm factor (\cite{hampton15,cohen18}). Note that for $d$ inputs, such a dimension is about $(w+1)^d$ for the tensor-product basis including at most the monomial of degree $w$.\\      
 
For models with moderate to large number of relevant inputs, Bayesian additive regression trees have been used for building emulators of such models using only the input and output observations (see \cite{horiguchi24} and references therein). Such approaches rely somehow on rule-ensemble learning approaches by constructing homogenous or local base learners (see \cite{friedman08} and references therein). Combining the model outputs and model derivatives can help for building emulators that account for both local and global properties of simulators. For instance, including derivatives in the Gaussian processes,  considered in \cite{solak02},  allows for improving such emulators. Emulators based on Taylor' series (see \cite{chkifa13}) combine  both the model outputs and derivative outputs with interesting theoretical results, such as dimension-free rates of convergence. However, concrete constructions of such emulators are not provided in that paper.\\      

 Note that the aforementioned emulators are part of the class of global approximations of functions. While global emulators are being able to be deployed for approximating models at any point of the entire domain, local or point-based emulators require building different emulators for the same model. 
Conceptually, the main issues related to such practical emulators are the truncation errors and the biases, due to the epistemic uncertainty. Indeed, none of the above emulators relies on exact and finite expansions of functions in general. Thus, additional assumptions about the order of such errors are necessary to derive the rates of convergence of emulators. For instance, decreasing eigenvalues of kernels are assumed in kernel-based leaning approaches (see e.g., \cite{tsybakov09}).\\ 
 So far, the recent derivative-based (Db) ANOVA provides exact expansions of smooth functions with different types of inputs' distributions using the model outputs, the first-order derivatives, second-order and up to the $d^{\mbox{th}}$-order cross partial derivatives (\cite{lamboni22}). It has been used in \cite{lamboni24stats} to derive the plug-in and crude emulators of models by replacing the unknown derivatives with their estimates. However, convergence analysis of such emulators are not known, and derivative-free methods are much convenient for applications in which the computations of cross-partial derivatives are too expansive or impossible (\cite{agarwal10,bach16,akhavan20}).   \\   
												     
 Therefore, for high-dimensional simulators for which all the inputs are important or not, it is worth investigating the development of their emulators based directly on available derivatives or derivative-free methods.  The contribution of this paper is threefold:     
\begin{itemize}   
\item  design adequate structures of emulators based on information gathered from global and derivative-based sensitivity analysis, such as unbiased orders of truncations and the selection of relevant ANOVA components (inputs and interactions); 
\item construct derivative-based or derivative-free global emulators that are easy to fit and compute, and can cope with every distribution of continuous input variables; 
\item examine the convergence analysis of our emulators with a particular focus on i) dimension-free upper-bounds of the biases and MSEs; ii) the parametric rates of convergence (i.e., $\mathsf{O}(N^{-1})$); and iii) the number of model runs needed to obtain the stability and accuracy of our estimations.        
\end{itemize}                     
         
In this paper, flexible emulators of complex models or simulators and approximations of functions are proposed using exact and finite expansions of functions, involving cross-partial derivatives, and known as Db-ANOVA. Section \ref{sec:deriaov} firstly deals with the general formulation of Db-ANOVA using different continuous distribution functions, and emulators of models with available derivatives. Such emulators reach the parametric rates of convergence and their MSEs do not depend on the dimensionality (i.e., $d$). Secondly, adequate structures of emulators are investigated using Sobol' indices and their upper bounds, as the components of such emulators are interpretable as main effects and interactions of a given order. The orders of unbiased truncations have been derived, leading to the possibility of selecting  the ANOVA components that are going to be included in our emulators.\\    
For non-highly smooth models, and for high-dimensional simulators for which the computations of derivatives are too expansive or impossible, new efficient emulators of simulators have been considered in  Section \ref{sec:dfexp0}, including their statistical properties. Firstly, we provide such emulators under the assumption of quasi-uniform distributions of inputs so as to i) obtain practical conditions for using such emulators, and ii) derive emulators that enjoy dimension-free MSEs for particular distributions of inputs. Secondly, such an assumption is removed to cope with every distribution of inputs. The proposed emulators have dimension-free biases and reach the parametric rate of convergence as well. Numerical illustrations (see Section \ref{sec:appli}) and an application to a heat diffusion PDE model (see Section \ref{sec:realapp}) have been considered to show the efficiency of our estimators, and we conclude this work in Section~\ref{sec:con}.  
                         
\section*{General notation} 
For an integer $d >0$, denote with $\bo{\X} := (X_1, \ldots, X_d )$ a random vector of $d$ independent and continuous variables with marginal cumulative distribution functions (CDFs) $F_j$  and probability density functions (PDFs) $\rho_j,\, j=1, \ldots, d$.\\  
For a non-empty subset $u \subseteq \{1, \ldots, d \}$, $|u|$ stands for its cardinality and $(\sim u) := \{1, \ldots, d \}\setminus u$.  Also, $\bo{\X}_u := (\X_j, \, \forall\, j \in u)$ denotes a subset of such variables, and the  partition $\bo{\X} = (\bo{\X}_u, \bo{\X}_{\sim u})$ holds.  
Finally, we use $\norme{\cdot}$ for the Euclidean norm, $\esp[\cdot]$ for the expectation operator and $\var[\cdot]$ for the variance operator.      
      
\section{New insight into derivative-based ANOVA and emulators} \label{sec:deriaov} 
 Given an integer $\NN >0$ and an open set $\Omega \subseteq \R^d$, consider a weak partial differentiable function $\M : \Omega \to \R^\NN$ (\cite{zemanian87,strichartz94}) and a subset  $v \subseteq \{1, \ldots, d\}$ with $|v|>0$. Denote with $\mathcal{D}^{|v|}  \M := \left(\prod_{k \in v} \frac{\partial }{\partial x_{k}} \right) \M$, the $|v|^{\text{th}}$ weak cross-partial derivatives of each component of $\M$ w.r.t. each $\x_{k}$ with $k \in v$ and $\mathcal{L}^2(\Omega) := \left\{\M : \Omega \to \R^\NN\, : \, \esp\left[\norme{\M(\bo{\X})}^2\right] <~+\infty \right\}$, the Hilbert space of functions. Consider the following Hilbert-Sobolev space:      
$$       
\displaystyle      
W^{d,2} := \left\{\M \in \mathcal{L}^2(\Omega) : \,   \mathcal{D}^{|v|}  \M \in  \mathcal{L}^2(\Omega); \quad \forall\; |v| \leq d \right\} \, .       
$$
 In what follows, assume that: \\                    
  
\begin{assump}[A1]
 $\bo{X}$ is a random vector of independent and continuous variables, supported on an open domain $\Omega$. 
\end{assump}          
\begin{assump}[A2]
 $\M(\cdot)$ is a deterministic function with $\M(\cdot) \in W^{d,2}$.   
\end{assump}          
   
\subsection{Full derivative-based emulators}  
Under (A2), every sufficiently smooth unction $\M(\cdot)$ admits the derivative-based ANOVA (Db-ANOVA) expansion (see \cite{lamboni22,lamboni24stats}), that is,   $\forall\, \bo{x} \in \Omega$, 
\begin{equation} \label{eq:deriGene}                 
\displaystyle                                      
\M(\bo{x})  = \esp_{\bo{\X}'} \left[ \M(\bo{\X}') \right]  +
  \sum_{\substack{v,\, v \subseteq \{1, \ldots, d\} \\ |v|>0}}  \esp_{\bo{\X}'} \left[
 \mathcal{D}^{|v|} \M\left(\bo{\X}' \right) \prod_{k \in v} \frac{G_{k}(\X_{k}') - \indic_{[\X_{k}' \geq x_{k}]}}{g_{k} (\X_{k}') }  \right] \, ,          
\end{equation}
where $\bo{\X}' :=\left(\X_1', \ldots, \X_d' \right)$ is a random vector of independent variables, having the CDFs $\X'_j \sim G_j$ and the PDFs $g_j := \frac{d G_j}{d x_j'}$. Evaluating $\M(\cdot)$ at the random vector $\bo{\X}$, and taking $G_j=F_j$ yield the unique and orthogonal  Db-ANOVA decomposition of $\M(\bo{\X})$, that is,                        
\begin{equation}     \label{eq:derivief}  
\displaystyle 
 \M(\bo{\X}) = \esp \left[\M(\bo{\X}')\right]  +  \sum_{\substack{v \subseteq \{1, \ldots, d\}\\|v|>0}} \M_v(\bo{\X}_v) ;     
\quad 
\M_v(\bo{\X}_v)  := \esp_{\bo{\X}'} \left[ \mathcal{D}^{|v|}\M(\bo{\X}') \prod_{k \in v} \frac{F_{k}(\X_{k}') - \indic_{[\X_{k}' \geq \X_{k}]}}{\rho_{k} (\X_{k}') }  \right] \, .      
\end{equation}   
When analytical cross-partial derivatives are available or the derivative data-sets are observed (see \cite{morris93,solak02} and the references therein),  we are able to derive emulators of complex models that are time demanding,  bearing in mind the method of moments. Indeed, given  
 a sample  $\left\{ \bo{\X}_i' \right\}_{i=1}^N := \left\{ \left(\X_{i,1}', \ldots, \X_{i, d}' \right) \right\}_{i=1}^N$ from $\bo{\X}'$, and the associated sample of (analytical or observed) derivatives outputs, that is, 
$$
\left( \left\{ \mathcal{D}^{|v|}\M(\bo{\X}_i')\right\}_{i=1}^N,\;  \forall\,  v \subseteq \{1, \ldots, d\} \right)  \, ,    
$$ 
the consistent (full) emulator or estimator of $\M$ at any sample point $\bo{x}$ of $\bo{\X}$ is        
\begin{equation}  \label{eq:emu}      
\displaystyle     
\widehat{\M_{N}}(\bo{x}) := \frac{1}{N} \sum_{i=1}^N \sum_{\substack{v,\,v \subseteq \{1,\ldots, d\}\\|v| \geq 0}}  
\mathcal{D}^{|v|}\M(\bo{\X}_i')  \prod_{k \in v} \frac{F_{k}(\X_{i,k}') - \indic_{[\X_{i,k}' \geq  x_{k}]}}{\rho_{k} (\X_{i,k}') } \, ,                
\end{equation}          
 with $\mathcal{D}^{|\emptyset|}\M = \M$ and $\prod_{k \in \emptyset} c_k :=~1$  for every real $c_k$ by convention.  We can check that $\widehat{\M_{N}}(\bo{x})$ is an unbiased estimator, and it reaches the parametric mean squared error (MSE) rate of convergence, that is,
$
\esp\left[\left(\widehat{\M_{N}}(\bo{x})- \M(\bo{x}) \right)^2 \right] = \mathcal{O}(N^{-1})  
$. This rate is dimension-free, provided that $\var\left[\sum_{\substack{v,\,v \subseteq \{1,\ldots, d\}\\|v| \geq 0}}  
\mathcal{D}^{|v|}\M(\bo{\X}_1')  \prod_{k \in v} \frac{F_{k}(\X_{1,k}') - \indic_{[\X_{1,k}' \geq  x_{k}]}}{\rho_{k} (\X_{1,k}') } \right] <~+\infty$.\\    
    
For complex models without cross-partial derivatives, optimal estimators of such derivatives (i.e., $\widehat{\mathcal{D}^{|v|}\M}$) have been used for constructing the plug-in consistent emulator of $\M$ in \cite{lamboni24stats}. Such an emulator is given by     
\begin{equation}  \label{eq:emu2}        
\displaystyle     
\widehat{\M_{N,p}}(\bo{x}) := \frac{1}{N} \sum_{i=1}^N \sum_{\substack{v,\,v \subseteq \{1,\ldots, d\}\\|v| \geq 0}}  
\widehat{\mathcal{D}^{|v|}\M}(\bo{\X}_i')  \prod_{k \in v} \frac{F_{k}(\X_{i,k}') - \indic_{[\X_{i,k}' \geq  x_{k}]}}{\rho_{k} (\X_{i,k}') } \, .   \nonumber                 
\end{equation}    
While the estimator $\widehat{\mathcal{D}^{|v|}\M}$ provided in \cite{lamboni24stats} has a dimension-free upper-bound of the bias, and reaches the parametric rate of convergence, its MSE increases with $d^{|v|}$, showing the necessary of using the number of model runs $N \propto  d^{|v|}$ to expect reducing significantly the MSE for higher-order cross-partial derivatives.  
   
\subsection{Adequate structures of emulators and truncations}  \label{sec:adstr}
In high-dimensional settings, it is common to expect reducing the dimensionality before building emulators. Truncations are common practice in polynomial approximation of functions (\cite{dongbin02,hampton15,wahba90,wong16}) and in ANOVA-practicing community (\cite{caflisch97,owen03,rabitz99,kucherenko09}), leading to truncated errors. Using the Db-ANOVA, controlling such errors  can be made possible according to the information gathered from global sensitivity analysis (\cite{kucherenko09,lamboni24stats}). Indeed, the variances of the terms in Db-ANOVA expansions of functions are exactly the main and interactions Sobol' induces up to a normalized constant when $G_j=F_j$. Thus, we are able to avoid non-relevant terms in our emulators according to the values of Sobol' induces, suggesting that adequate truncations will not have any impact on the MSEs and the above parametric rate of convergence. For the sake of simplicity, $\NN=1$ is considered in what follows.     

\begin{defi}
Consider an integer $d_0 \in \{1, \ldots, d\}$ and the full Db-ANOVA given by (\ref{eq:derivief}). The  truncated Db-ANOVA of $\M$ (in the superpose sense (\cite{caflisch97,owen03,rabitz99,kuo10})) of order $d_0$ is given by 
$$
\displaystyle  
 \M_{T,d_0}(\bo{\X}) := \esp \left[\M(\bo{\X}')\right]  +  \sum_{\substack{v \subseteq \{1, \ldots, d\}\\ 0< |v| \leq d_0}} \M_v(\bo{\X}_v)  \, . 
$$   
\end{defi}  
While $\M_{T,d_0}$ is an approximation of $\M$ in general, the equality holds for some  functions. Given an integer $\alpha \geq 0$,  consider the space of functions   
$$  
\displaystyle           
\mathcal{L}_{\alpha, 0} := \left\{ 
\M : \R^d \to \R^\NN \, : \,
\left|\M(\bo{x}) -
\sum_{\substack{w \subseteq \{1, \ldots, d\} \\ |w| \leq \alpha}} \esp_{\bo{\X}'}\left[ \mathcal{D}^{|w|}\M(\bo{\X}') \prod_{k \in w} \frac{G_k \left(\X_k' \right) -\indic_{\X_k' \geq x_k }}{g_k\left(\X_k' \right)} \right] \,  \right| =0 \right\} \, .    
$$      
We can see that $\mathcal{L}_{\alpha, 0}$ is a class of functions having at most $\alpha$-order  interactions. Also, $\mathcal{L}_{\alpha, 0}$ contains the class of functions given by 
$$
 \mathcal{D}^{|\{v, \jmath_0\}|}\M = 0, \;  \forall\, v \subseteq \{1, \ldots, d\}, \, 
\jmath_0 \in (\sim v),  \, |v|= \alpha   \, ,   
$$       
It is then clear that we have $\M(\bo{\X}) = \M_{T,d_0}(\bo{\X})$ when $\M \in \mathcal{L}_{\alpha=d_0, 0}$.
\begin{defi}   
A truncation of $\M$ given by $\M_{T,d_0}$ is said to be an unbiased one whenever $\M(\bo{\X}) = \M_{T,d_0}(\bo{\X})$. 
\end{defi} 
Thus, $\mathcal{L}_{\alpha, 0}$ is the class of unbiased  truncations of $\M$ up to the order $d_0=\alpha$.\\  

Formally, based on available derivatives, we are able to derive unbiased  truncations for some classes of functions. Consider the following Db-expressions of Sobol' indices of the input variables $\X_j,\, j=1, \ldots, d$ and their upper-bounds (see \cite{lamboni22}):   
$$
S_{j} := \frac{1}{\var\left[\M(\bo{\X}) \right]} \esp\left[\frac{\partial\M}{\partial x_j} (\bo{\X})  \frac{\partial\M}{\partial x_j} (\bo{\X}')\frac{ F_j\left(\min(\X_j, \X_j')\right) - F_j(\X_j) F_j(\X_j')}{\rho_j(\X_j) \rho_j(\X_j')} \right]  \, ,    
$$ 
$$
 S_{T_j} :=  \frac{1}{\var\left[\M(\bo{\X}) \right]} \esp\left[\frac{\partial\M}{\partial x_j} (\bo{\X})  \frac{\partial\M}{\partial x_j} (\X_j', \bo{\X}_{\sim j})\frac{ F_j\left(\min(\X_j, \X_j')\right) - F_j(\X_j) F_j(\X_j')}{\rho_j(\X_j) \rho_j(\X_j')} \right]  \, ,  
$$
and   
$$
S_j \leq 
S_{T_j} \leq  U\!B_j := \frac{1}{2 \var\left[\M(\bo{\X}) \right]} \esp\left[\left(\frac{\partial\M}{\partial x_j} (\bo{\X}) \right)^2 \frac{ F_j(\X_j) \left(1- F_j(\X_j)\right)}{\left[\rho_j(\X_j)\right]^2} \right] \, . 
$$
Note that the computations of $S_{j}$ and $ U\!B_j$ are straightforward using given first-order derivatives (see \cite{lamboni24stats}), while those of $S_{T_j}$ require using i) non-parametric methods for given derivatives or ii) derivatives of specific input values. Using such indices, adequate structures of $\M(\cdot)$ can be constructed. For instance, it is known that when $\sum_{j=1}^d S_{j} = 1$,  $\M$ is an additive function of the form $\M(\bo{x}) =\sum_{j=1}^d h_{j}(x_j)$ with $h_j$ a real valued function. Thus, a truncation (in the superpose sense) of order $d_0=1$ is an unbiased truncation. Other values of $d_0$ are given in Proposition \ref{prop:sedstrc}.  
  
\begin{prop}   \label{prop:sedstrc}
Consider the main and total indices given by $S_j, \, S_{T_j}$ with $j=1, \ldots, d$. Then,   
\begin{itemize}
\item $\sum_{j=1}^d ( S_j + S_{T_j})  = 2$ implies that using $d_0=2$ leads to unbiased truncations;   
\item $\sum_{j=1}^d ( 2S_j + S_{T_j})  \in ]2,\, 3]$ implies that  using $d_0=3$ leads to unbiased truncations;    
\item  if there exists an integer $\alpha$ such that  $(\alpha-1) \sum_{j=1}^d S_{T_j} +  \sum_{j=1}^d  S_{j} \geq \alpha$ and  $\sum_{j=1}^d  S_{T_j} + (\alpha-1) \sum_{j=1}^d  S_{j} \leq \alpha$, then $d_0=\alpha$ leads to unbiased truncations.    
\end{itemize}     
\end{prop}       
\begin{preuve}
See Appendix \ref{app:prop:sedstrc}.       
\hfill $\square$
\end{preuve}  

In general, if $
 \mathcal{D}^{|\{v, \jmath_0\}|}\M = 0, \;  \forall\, v \subseteq \{1, \ldots, d\}, \, 
\jmath_0 \in (\sim v),  \, |v|= \alpha         
$,           
then  $d_0=\alpha$ leads to unbiased truncations. Often, low-order derivatives are available or can be efficiently computed using fewer model runs, leading to truncated emulators in the superpose sense. It is worth noting that our emulators still enjoy the above parametric rate of convergence for unbiased truncations. In the presence of the truncated errors, it is usually difficult to derive the rates of convergence without additional assumptions on the order of such errors.\\ 
 
In addition to such truncations in the superpose sense,  screening measures allow for quickly identifying non-relevant inputs (i.e., $U\!B_j \approx 0$), leading to possible dimension reductions. For instance, we can see that    
\begin{itemize}  
\item  $S_{j} = S_{T_j}$ or $S_{j} \approx U\!B_j $ implies removing all cross-partial derivatives or interactions involving $\X_j$;  
\item   $S_{j} = 0 $ and $U\!B_j \geq S_{T_j} \gg 0, \, \forall\, j \in \{1, \ldots, d \}$, suggest removing the first-order terms, corresponding to $d_0=1$;    
\item $1 -(S_i + S_{T_j}) \leq S_{T_{\sim \{i,j\}}} = 0$ or equivalently $S_i + S_{T_j}=1$ and $S_{T_k} =0,\, \forall\, k\notin \{i, j\}$ implies keeping only $\X_i$ and $\X_j$.    
\end{itemize}                
 For non-highly smooth functions and the models for which the computations of derivatives are too expansive or impossible, derivative-free methods combined with unbiased truncations remain an interesting framework.   
        
\section{Derivative-free emulators of models} \label{sec:dfexp0} 
 This section aims at developing emulators of models, even when all the inputs are important according to screening measures.
           
\subsection{Stochastic surrogates of functions using Db-ANOVA} \label{sec:dfexp} 
Consider integers $L >0, \, q >0$; $\beta_{\ell} \in \R$ with $\ell=1, \ldots, L$; $\boldsymbol{\hh} := (\hh_1, \ldots, \hh_d) \in \R^d_+$, and  denote with $\bo{V} :=(V_1, \ldots, V_d)$ a $d$-dimensional random vector of independent variables  satisfying:  $\forall \, j \in \{1, \ldots, d\} $, 
$$ 
 \esp\left[V_j\right] =0; \qquad \esp\left[\left(V_{j} \right)^{2}\right] =\sigma^2;   \qquad \esp\left[\left(V_{j} \right)^{2q+1}\right] =0; \qquad \esp\left[\left(V_{j} \right)^{2q}\right] < +\infty \, .  
$$  
Any random variable which is symmetric about zero is an instance of $V_j$s. Also, denote $\beta_\ell \boldsymbol{\hh}\bo{V}  := (\beta_{\ell}\hh_1 V_1; \ldots, \beta_{\ell} \hh_d V_d)$.          
For concise reporting of the results, the elementary symmetric polynomials (ESP) are going to be used (see e.g.,  \cite{alatawi22,ahmed23}).       
\begin{defi}  
Given $u \subseteq \{1, \ldots, d\}$ with $|u|>0$ and $\bo{r}_u := (r_k, \forall\, k \in u) \in \R^{|u|}$, the $p^{\mbox{th}}$ ESP of $\bo{r}_u$ is defined as follows: 
$$  
\bo{e}_p^{(u)}(\bo{r}_u) := \left\{ \begin{array}{cl}
0 & \mbox{if} \; \;  p> |u| \: \mbox{or} \; p<0 \\
1 & \mbox{if} \; \;  p =0 \\
\sum_{\substack{w \subseteq u \\ |w|=p}} \prod_{k \in w} r_k & \mbox{if} \; \;   p =1, \ldots, |u| \\    
\end{array} \right. \, .         
$$          
\end{defi}    
Note that $\bo{r} := \bo{r}_{\{1, \ldots, d\}} = (r_1, \ldots, r_d) \in \R^d$, $\bo{e}_p^{(1:d)}\left( \bo{r}\right) :=\bo{e}_p^{(\{1, \ldots, d\})}\left( \bo{r}\right)$.  Also, given $\X_j' \sim G_j, \; \forall\,  j \in \{1, \ldots, d\}$,  define      
$$
R_k\left(x_{k}, \X_{k}', V_k \right) := \frac{G_{k}(\X_{k}') - \indic_{[\X_{k}' \geq x_{k}]}}{g_{k} (\X_{k}')\, \hh_k  \sigma^2 } V_k \; , \qquad k=1, \ldots, d \, ;   
$$      
$$
\bo{R}_u\left(\bo{x}_u, \bo{\X}'_u, \bo{V}_u \right) := \left( R_k\left(x_{k}, \X_{k}', V_k \right), \forall\, k \in u \right); \quad \bo{R}\left(\bo{x}, \bo{\X}', \bo{V} \right) :=   \bo{R}_{\{1, \ldots, d\}}\left(\bo{x}, \bo{\X}', \bo{V} \right)  \, .    
$$  

Without loss of generality, we are going to focus on the modified output, that is, $\M^c(\bo{x}) := \M(\bo{x}) -\esp\left[\M(\bo{\X}') \right]$. Based on the above framework, Theorem \ref{theo:sfall} provides a new approximation of every function or surrogate of a deterministic simulator.                
     
\begin{theorem} \label{theo:sfall}   
Consider distinct $\beta_\ell$s. If $\M$ is smooth enough and (A1) hold, then there exists  $\alpha_{d} \in \{1, \ldots, L\}$ and real coefficients $C_{1}^{(p)}, \ldots, C_{L}^{(p)}, \, p=1, \ldots, d$ such that  
\begin{equation}  \label{eq:ceout}     
\displaystyle  
\M^c(\bo{x})  = \sum_{\ell=1}^{L} \sum_{p=1}^{d}  C_{\ell}^{(p)}  
  \esp \left[ \M\left(\bo{\X}' + \beta_\ell \boldsymbol{\hh}\bo{V} \right) \bo{e}_p^{(1:d)}\left( \bo{R}\left(\bo{x}, \bo{\X}', \bo{V} \right) \right) \right] +   \mathcal{O}\left( \norme{\bo{\hh}}^{2\alpha_d} \right)   \, .
\end{equation}         
\end{theorem}                     
\begin{proof}    
See Appendix \ref{app:theo:sfall}.      
\end{proof}      
  
The setting $L=1, \beta_1 =1, \, C_{1}^{(p)}=1$  with $p=1, \ldots, d$ is going to provide an approximation of order $\mathcal{O}\left( \norme{\bo{\hh}}^{2} \right)$. Equivalently, the same order is obtained  using  the constraints  
\begin{equation}    
\left\{  
\begin{array}{ll}         
\sum_{\ell=1}^{L} C_{\ell}^{(p)} \beta_{\ell}^{r}= \delta_{p,r}; \quad  r =0, \ldots, L-1,  & 
\mbox{if} \; p \leq  L-1 \\  
\sum_{\ell=1}^{L} C_{\ell}^{(p)} \beta_{\ell}^{r}= \delta_{p,r}; \quad  r =0, \ldots, L-2, p,   & \mbox{otherwise} \\
\end{array}        
\right.  \, ,          \nonumber                
\end{equation}  
 with $p=1, \ldots, d$ and  $L \leq d+1$. In the same sense, taking     
$  
\sum_{\ell=1}^{L} C_{\ell}^{(p)} \beta_{\ell}^{r}= \delta_{p,r}; \;  r =p, \ldots, L+p-1$ with $p=1, \ldots, d$ yields an approximation of order $\mathcal{O}\left( \norme{\bo{\hh}}^{2L} \right)$. Such constraints implicitly define the coefficients $C_{1}^{(p)}, \ldots, C_{L}^{(p)}, \, p=1, \ldots, d$, and  they rely on the (generalized) Vandermonde matrices. Distinct values of $\beta_\ell$s (i.e., $ \beta_{\ell_1} \neq \beta_{\ell_2} $) ensure the existence and uniqueness of such coefficients, as such matrices are invertible (see \cite{rawashdeh19,ahmed23}).   \\   

To improve the approximations of lower-order terms (i.e., lower values of $p$ ) in Equation (\ref{eq:ceout}), we are given an integer $r^* \in \{0, \ldots, d-1\}$ with $r^* \leq L-2$ and consider the following constraints:  
\begin{equation}  \label{eq:consttype1}         
\left\{  
\begin{array}{ll}         
\sum_{\ell=1}^{L} C_{\ell}^{(p)} \beta_{\ell}^{r}= \delta_{p,r}; \quad  r =0, \ldots, r^*,  p+2\lambda_p+2, \ldots, p+2\lambda_p +2(L-1-r^*),  &    
\mbox{if} \; p \leq r^* \\     
\sum_{\ell=1}^{L} C_{\ell}^{(p)} \beta_{\ell}^{r}= \delta_{p,r}; \quad  r =0, \ldots, r^*, p, p+2 \ldots, p+ 2(L-r^*-2),  & \mbox{otherwise} \\                       
\end{array}        
\right.  \, ,                 
\end{equation}     
where $\lambda_p := \left[\frac{r^*-p}{2} \right]$ stands for the largest integer that is less than $\frac{r^*-p}{2} $. The above choice of coefficients requires $L\geq r^*+2$, and $L^* := r^*+2$ is the minimum number of model runs used for deriving surrogates of functions. Such coefficients are much suitable for truncated surrogates. For instance, the truncated surrogate of order $d_0 \leq d$ (in the superposition sense) is given by  
$$
\displaystyle  
\widetilde{\M_{T, d_0}^c}(\bo{x})  := \sum_{\ell=1}^{L} \sum_{p=1}^{d_0}  C_{\ell}^{(p)}  
  \esp \left[ \M\left(\bo{\X}' + \beta_\ell \boldsymbol{\hh}\bo{V} \right) \bo{e}_p^{(1:d)}\left( \bo{R}\left(\bo{x}, \bo{\X}', \bo{V} \right) \right) \right]   \, .    
$$
Note that in this case, one must require $r^* \in \{0, \ldots, d_0-1\}$, and we will see that $r^*=d_0-1$ is the best choice to improve the MSEs. \\  

Likewise, when $\bo{\X}_{u_I}$ with $u_I \subset \{1, \ldots, d\}$ is the vector of the most influential input variables according to variance-based sensitivity analysis (see Section \ref{sec:adstr}), the following truncated surrogate should be considered:  
$$
\widetilde{\M_{u_I}^c}\left(\bo{x}_{u_I} \right) := \sum_{\ell=1}^{L} \sum_{p=1}^{|u_I|}  C_{\ell}^{(p)}  
  \esp \left[ \M\left(\bo{\X}' + \beta_\ell \boldsymbol{\hh}\bo{V} \right) \bo{e}_p^{(u_I)}\left( \bo{R}_{u_I}\left(\bo{x}_{u_I}, \bo{\X}'_{u_I}, \bo{V}_{u_I} \right) \right) \right]\, .                
$$        
 
Based on Equation (\ref{eq:ceout}), the method of moments allows for deriving the emulator of any simulator or the estimator of any function. To that end, we are given two independent samples of size $N$, that is, $\left\{\bo{\X}_i' \right\}_{i=1}^N :=\left\{ \left(\X_{i,1}',  \ldots, \X_{i,d}'\right) \right\}_{i=1}^N$ from $\bo{\X}'$ and  $\left\{\bo{V}_i \right\}_{i=1}^N := \left\{ \left(V_{i,1},  \ldots, V_{i,d}\right) \right\}_{i=1}^N$ from $\bo{V}$. The full and consistent emulator is given  by 
\begin{equation}  \label{eq:emuappdfful}   
\displaystyle       
\widehat{\M^{c}_N}(\bo{x}) := \frac{1}{N} \sum_{i=1}^N  
\sum_{\ell=1}^{L} \sum_{p=1}^{d}  C_{\ell}^{(p)}  
   \M\left(\bo{\X}'_i + \beta_\ell \boldsymbol{\hh}\bo{V}_i \right) \bo{e}_p^{(1:d)}\left( \bo{R}\left(\bo{x}, \bo{\X}'_i, \bo{V_i} \right) \right) \, .     
\end{equation}    
The derivations of the truncated emulators (i.e., $\widehat{\widetilde{\M_{T, d_0}^c}}$ and $\widehat{\widetilde{\M_{u_I}^c}}$) are straightforward. All these emulators rely on $NL$ model runs with the possibility $L\ll d$. This property is useful for high-dimensional simulators.    
       
\subsection{Statistical properties of our emulators} \label{sec:statprop} 
While the emulator $\widehat{\M^{c}_N}(\bo{x})$ does not rely on the model derivatives, structural and technical assumptions on $\M$ are needed for deriving the biases of this emulator, such as the  H\"older  space of functions. Given $\Vec{\boldsymbol{\imath}} :=(i_1, \ldots, i_d)  \in \N^d$, denote    
$\mathcal{D}^{(\Vec{\boldsymbol{\imath}})}\M := \left(\prod_{k=1}^d  \frac{\partial^{i_k} }{\partial x_{k}} \right)\M$, $(\bo{x})^{\Vec{\boldsymbol{\imath}}} := \prod_{k=1}^d x_k^{i_k}$, $\Vec{\boldsymbol{\imath}}! = i_1! \ldots i_d!$ and $||\Vec{\boldsymbol{\imath}}||_1 =i_1+\ldots i_d$. Given $\alpha \geq 0$, the H\"older  space of $\alpha$-smooth functions is given by $\forall \, \bo{x}, \bo{y} \in \R^d$,    
$$        
\displaystyle           
\mathcal{H}_\alpha := \left\{ 
\M : \R^d \to \R\, : \,
\left|\M(\bo{x}) - \sum_{||\Vec{\boldsymbol{\imath}}||_1=0}^{\alpha-1} 
\frac{\mathcal{D}^{(\Vec{\boldsymbol{\imath}})}\M(\bo{y})}{\Vec{\boldsymbol{\imath}}!} \left(\bo{x}-\bo{y} \right)^{\Vec{\boldsymbol{\imath}}} \right|  
  \leq M_\alpha \norme{\bo{x}- \bo{y}}^\alpha \right\} \, ,     
$$            
with $M_\alpha>0$ and $\mathcal{D}^{(\Vec{\boldsymbol{\imath}})}\M(\bo{y})$ a (weak) cross-partial derivative.\\

Also, given $B_{\alpha} \geq 0$ and CDFs $G_j$s,  define the following space of functions:  
$$
\displaystyle           
\mathcal{L}_{\alpha, B_{\alpha}} := \left\{ 
\M : \R^d \to \R\, : \,
\left|\M(\bo{x}) -
\sum_{\substack{w \subseteq \{1, \ldots, d\} \\ |w| \leq \alpha}} \esp_{\bo{\X}'}\left[ \mathcal{D}^{|w|}\M(\bo{\X}') \prod_{k \in w} \frac{G_k \left(\X_k' \right) -\indic_{\X_k' \geq x_k }}{g_k\left(\X_k' \right)} \right] \,  \right| \leq  B_{\alpha} \right\} \, .   
$$      
We can see that $\mathcal{L}_{\alpha, B_{\alpha}}$ contains constants; $\mathcal{L}_{\alpha, 0}$ is a class of functions having at most $\alpha$-order of interactions;  and $\mathcal{L}_{d, 0}$ is a class of all smooth functions, as $B_{d}=0$.  Lemma \ref{lem:trucbias} provides the links between both spaces. To that end, consider $M_{|w|}' := \norminf{\mathcal{D}^{|w|}\M}$ for all $w \subseteq \{1, \ldots, d\} $ with $\norminf{\cdot}$ the infinity norm. 

\begin{assump}[A3]   
 $g_j \geq \rho_{\min} >0$ for any $j \in \{1, \ldots, d\}$ \, .       
\end{assump} 

Assumption (A3) aims at covering the class of quasi-uniform distributions, and other distributions for which the event $g_j \geq \rho_{\min}$ occurs with a high probability. It is the case of most unbounded distributions. 
    
\begin{lemma} \label{lem:trucbias} 
Consider $0 <d_0 \leq d$, and assume that $\M \in \mathcal{H}_{\alpha}$ with $\alpha \in \{0, d \}$ and (A1), (A3) hold. Then, there exists $\gamma_0 >0$ such that $\M \in \mathcal{L}_{d_0, D_{d_0, \rho_{\min}}}$, with
\begin{equation} \label{eq:disc} 
D_{d_0, \rho_{\min}} := 
2 \gamma_0 M_{0}' \left[\frac{1}{2  \rho_{\min}} \left(\frac{M_{d}'}{M_{0}'} \right)^{1/d} +1 \right]^{d_0}
\left\{ \left[\frac{1}{2  \rho_{\min}} \left(\frac{M_{d}'}{M_{0}'} \right)^{1/d} +1 \right]^{d-d_0} -1 \right\} \, . 
\nonumber  
\end{equation}     
\end{lemma}   
\begin{preuve}
See Appendix \ref{app:lem:trucbias}.   
\hfill $\square$
\end{preuve}     

Remark that Lemma \ref{lem:trucbias} also provides the upper-bound of the remaining terms when approximating $\M(\bo{x})$ by the truncated function 
$$
\M_{d_0}(\bo{x}) := \sum_{\substack{w \subseteq \{1, \ldots, d\} \\ |w| \leq d_0}} \esp_{\bo{\X}'}\left[ \mathcal{D}^{|w|}\M(\bo{\X}') \prod_{k \in w} \frac{G_k \left(\X_k' \right) -\indic_{\X_k' \geq x_k }}{g_k\left(\X_k' \right)} \right] \, .
$$ 
 When $ \frac{1}{2  \rho_{\min}} \left(\frac{M_{d}'}{M_{0}'} \right)^{1/d} \to~0$, $D_{d_0, \rho_{\min}}$ is equivalent to   
$$
D_{d_0, \rho_{\min}} \equiv  \frac{(d-d_0) \gamma_0 \left(M_{0}'\right)^{(d-1)/d} \left(M_{d}'\right)^{1/d} }{\rho_{\min}}\left[\frac{d_0}{2  \rho_{\min}} \left(\frac{M_{d}'}{M_{0}'} \right)^{1/d} +1 \right] \, . 
$$  

\subsubsection{Biases of the proposed emulators}
To derive the bias of $\widehat{\M^c_{N, d_0}}$ (i.e., estimator of $\M^c_{d_0}$) in Theorem \ref{theo:biasfc} using the aforementioned  spaces of functions, denote with $\bo{Z} :=(Z_1, \ldots, Z_d)$ a $d$-dimensional random vector of independent variables that are centered about zero and standardized (i.e., $\esp[Z_k^2]=1$, $k=1, \ldots, d$), and $\mathcal{R}_c$ the set of such random vectors. For any $r \in \N$ and $w \subseteq \{1, \ldots, d\}$,  define        
$$     
\Gamma_{r} := \sum_{\ell=1}^{L}  \left| C_\ell^{(|w|)} \beta_\ell^{r}   \right|;
\qquad \qquad 
K_{w,L} := \inf_{\bo{Z} \in \mathcal{R}_c} \esp\left[ \norme{\bo{Z}^2}^{L} \prod_{k \in w} Z_k^2 \right] 
\Gamma_{|w|+2L}  \, ;           
$$       
$$
L'_w := \left(\left[\frac{r^*-|w|}{2} \right] + L-r^* \right) \indic_{|w|\leq r^*} +  (L-r^*-1) \indic_{|w|> r^*}
\, . 
$$        
              
\begin{theorem}   \label{theo:biasfc}    
 Assume $\M \in \mathcal{H}_{\alpha}$ with $\alpha \in \left\{0, \max(d, d_0+2(L-r^*-1)) \right\}$ and (A1), (A3) hold. Then, we have 
\begin{equation}
\left|\esp\left[\widehat{\M^c_{N, d_0}}(\bo{x}) \right] -\M^c(\bo{x}) \right| \leq  
\sum_{\substack{w \subseteq \{1, \ldots, d\} \\ 0< |w| \leq  d_0}}   \sigma^{2L'_w}  \norme{\boldsymbol{\hh}^2}^{L'_w }  M_{|w|+2L'_w } K_{w,L'_w }  \left(\frac{1}{2  \rho_{\min}} \right)^{|w|} + D_{d_0, \rho_{\min}} \, .  
\end{equation}      
 Moreover, if  $V_k \sim \mathcal{U}(-\xi, \xi)$ with $\xi>0$ and $k=1, \ldots, d$, then 
\begin{equation}
\left|\esp\left[\widehat{\M^c_{N, d_0}}(\bo{x}) \right] -\M^c(\bo{x}) \right|  \leq  
\sum_{\substack{w \subseteq \{1, \ldots, d\} \\ 0< |w| \leq  d_0}}  \xi^{2L'_w }  ||\boldsymbol{\hh}^2||_1^{L'_w}  
 M_{|w|+2L'_w }  \Gamma_{|w|+2L'_w}  \left(\frac{1}{2  \rho_{\min}} \right)^{|w|} + D_{d_0, \rho_{\min}}  \, .     
\end{equation}     
\end{theorem}      
\begin{preuve}
See Appendix \ref{app:theo:biasfc}.  
\hfill $\square$ 
\end{preuve}

Using the fact $\hh_k \to 0$, the results provided in Theorem \ref{theo:biasfc} have simple upper-bounds (see Corollary~\ref{coro:biasfc1}). To provide such results, consider 
$$
 K_{1,r^*, d_0}^{\max} := \max_{\substack{w \subseteq \{1, \ldots, d\} \\ r^* < |w| \leq  d_0}} \left\{ K_{w, (L-r^*-1)}  M_{|w|+ 2(L-r^*-1) } \right\} \, ; 
$$
$$
K_{2, r^*, d_0}^{\max} := \max_{\substack{w \subseteq \{1, \ldots, d\} \\ r^* < |w| \leq  d_0}} \left\{  M_{|w|+ 2(L-r^*-1)} \Gamma_{|w|+2(L-r^*-1)} \right\} \, ;       
$$      
$$
K_{1, \rho_{\min}, r^*} := \left[ 2  \rho_{\min} \left(\frac{d}{2  \rho_{\min}}\right)^{r^*+1} \frac{ \left(\frac{d}{2  \rho_{\min}}\right)^{d_0-r^*} - 1}{d-2  \rho_{\min}} \right] \indic_{r^* < d_0-1} +  \binom{d}{d_0} \left(\frac{1}{2  \rho_{\min}} \right)^{d_0} \indic_{r^*= d_0-1} \, .   
$$           
\begin{corollary}  \label{coro:biasfc1}  
  Assume $\M \in \mathcal{H}_{\alpha}$ with $\alpha \in \left\{0, \max(d, d_0+2(L-r^*-1)) \right\}$ and (A1), (A3) hold. If $\hh_k \to 0$, then    
\begin{equation}  
\left|\esp\left[\widehat{\M^c_{N, d_0}}(\bo{x}) \right] -\M^c(\bo{x}) \right| \leq    \norme{\boldsymbol{\hh}^2}^{L-r^*-1}
  \sigma^{2(L-r^*-1)}  
 K_{1,r^*, d_0}^{\max} K_{1, \rho_{\min}, r^*}   +  D_{d_0, \rho_{\min}}    +   \mathcal{O}\left(\norme{\boldsymbol{\hh}^2}^{L-r^*-1} \right) \, .  
\end{equation}                       
 Moreover, if  $V_k \sim \mathcal{U}(-\xi, \xi)$ with $\xi>0$ and $k=1, \ldots, d$, then 
\begin{equation}
\left|\esp\left[\widehat{\M^c_{N, d_0}}(\bo{x}) \right] -\M^c(\bo{x}) \right| \leq  ||\boldsymbol{\hh}^2||_1^{L-r^*-1}  \xi^{2(L-r^*-1)}  K_{2, r^*, d_0}^{\max} K_{1, \rho_{\min}, r^*}   + D_{d_0, \rho_{\min}}  + \mathcal{O}\left( ||\boldsymbol{\hh}^2||_1^{L-r^*-1} \right)   \, .     
\end{equation}        
\end{corollary}        
\begin{preuve}
See Appendix \ref{app:coro:biasfc1}. 
\hfill $\square$     
\end{preuve} 

Using the above results, the bias of the full  emulator of $\M^c$ is straightforward by taking $d_0=d$ and knowing that $D_{d, \rho_{\min}}=0$. Moreover, Corollary \ref{coro:biasfc2} provides the bias of such an emulator under different structural assumptions on $\M$ so as to cope with many functions. To that end, define 
$$  
K_{1:d} := \inf_{\bo{Z} \in \mathcal{R}_c} \esp\left[ \norme{\bo{Z}} \prod_{k=1}^d Z_k^2 \right] \, .
$$      
         
\begin{corollary}  \label{coro:biasfc2} 
Let $d_0=d$; $r^* \leq d-1$ and $L=r^*+2$.  Assume $\M \in \mathcal{H}_{\alpha}$ with $\alpha \in \left\{0, d+1 \right\}$ and (A1), (A3) hold. If $\hh_k \to 0$, then 
\begin{equation}  
\left|\esp\left[\widehat{\M^c_{N}}(\bo{x}) \right] -\M^c(\bo{x}) \right| \leq    \sigma  \norme{\boldsymbol{\hh}}  M_{d+1} K_{1:d}  \, \Gamma_{d+1} \,  \prod_{k=1}^d  \esp \left[\left| E_k \right| \right]  +  \mathcal{O}\left(\norme{\boldsymbol{\hh}} \right) \, . 
\end{equation}                   
 Moreover, if  $V_k \sim \mathcal{U}(-\xi, \xi)$ with $\xi>0$ and $k=1, \ldots, d$, then 
\begin{equation} \label{eq:upful}
\left| \esp\left[\widehat{\M^c_{N}}(\bo{x}) \right] -\M^c(\bo{x}) \right| \leq    \xi ||\boldsymbol{\hh}||_1   M_{d+1} \, \Gamma_{d+1} \, 
\prod_{k=1}^d   \esp \left[ \left| E_k \right| \right]  + \mathcal{O}\left(||\boldsymbol{\hh}||_1 \right)   \, .     
\end{equation}             
\end{corollary}        
\begin{preuve}
See Appendix \ref{app:coro:biasfc2}.  
\hfill $\square$  
\end{preuve}  

In view of Corollaries \ref{coro:biasfc1}-\ref{coro:biasfc2},  Equation (\ref{eq:upful}) can lead to a dimension-free upper-bound of the bias. Indeed, using the uniform bandwidth $\hh_k=\hh$ and 
$$
\xi \leq \frac{1}{d M_{d+1} \, \Gamma_{d+1} \prod_{k=1}^d   \esp \left[ \left| E_k \right| \right]} \, , 
$$
we can see that the upper-bound of the bias is $\left| \esp\left[\widehat{\M^c_{N}}(\bo{x}) \right] -\M^c(\bo{x}) \right| \leq   \hh$. 
Furthermore, it is worth noting that $\left|\esp\left[\widehat{\M^c_{N, d_0}}(\bo{x}) \right] -\M^c(\bo{x}) \right| \leq \hh^{L-r^*-1}$ for any function $\M \in \mathcal{L}_{d_0,0}$.   
         
\subsubsection{Mean squared errors}
We start this section with the variance of the proposed emulators, followed by their mean squared errors and different rates of convergence. For the variance, define 
$$  
\digamma_{r^*,d_0}^{\max} := \max_{\substack{w \subseteq \{1, \ldots, d\} \\ |w|\leq d_0}} \left\{ M_{\min(r^*, |w|-1) +1} \Gamma_{\min(r^*, |w|-1) +1}\right\};
\qquad \quad 
\digamma_{d_0}^{\max} := \max_{\substack{w \subseteq \{1, \ldots, d\} \\ |w|\leq d_0}} \left\{ M_{|w|} \Gamma_{|w|}\right\}  \, .       
$$ 
    
\begin{theorem}   \label{theo:rtestfc}    
Consider the coefficients given by Equation (\ref{eq:consttype1}), and assume $\M \in \mathcal{H}_{\alpha}$ with $\alpha \in \left\{0, \max(d, d_0+2(L-r^*-1)) \right\}$ and (A1)-(A3) hold. Then, 
$$
\var\left[\widehat{\M^c_{N, d_0}}(\bo{x}) \right]    \leq  \frac{\left(\digamma_{r^*,d_0}^{\max} \right)^2}{N}
\sum_{\substack{w \subseteq \{1, \ldots, d\} \\ 0< |w| \leq d_0}}
 \prod_{k \in w} \left( \frac{\esp\left[E_k^2\right]}{\hh_k^2 \sigma^4} \esp\left[ V_1^2 \norme{\boldsymbol{\hh}\bo{V} }^{2\frac{\min(r^*, |w|-1) +1}{|w|}}\right] \right)   \, . 
$$
 Moreover, if  $r^*=d_0-1$, $\hh_k=\hh$ and $Z_k =V_k/\sigma$, then 
$$
\var\left[\widehat{\M^c_{N, d_0}}(\bo{x}) \right]    \leq \frac{d \left(\digamma_{d_0}^{\max} \right)^2 \esp \left[Z_1^2 \norme{\bo{Z} }^{2} \right]}{N\left(d\esp \left[Z_1^2 \norme{\bo{Z} }^{2} \right] - 3 \rho_{\min}^2\right)}  \left[\left( \frac{d\esp \left[Z_1^2 \norme{\bo{Z} }^{2} \right]}{3 \rho_{\min}^2}   \right)^{d_0}-1 \right] \, . 
$$     
\end{theorem}       
\begin{preuve}
See Appendix \ref{app:theo:rtestfc}.        
\hfill $\square$
\end{preuve}    

It turns out that the upper-bounds of the variance in Theorem \ref{theo:rtestfc} do not depend on the uniform bandwidths when $r^*=d_0-1$, leading to the derivations of the parametric MSEs of $\widehat{\M^c_{N, d_0}}(\bo{x})$ and $\widehat{\M^c_{N}}(\bo{x})$. To that end, consider the upper-bound of the above variance, that is,      
$$   
\bar{\Sigma}_{d_0} := \frac{d \left(\digamma_{d_0}^{\max} \right)^2 \esp \left[Z_1^2 \norme{\bo{Z} }^{2} \right]}{N\left(d\esp \left[Z_1^2 \norme{\bo{Z} }^{2} \right] - 3 \rho_{\min}^2\right)}  \left[\left( d\frac{\esp \left[Z_1^2 \norme{\bo{Z} }^{2} \right]}{3 \rho_{\min}^2}   \right)^{d_0}-1 \right]  \, .     
$$    
  
\begin{rem}  \label{rem:chvjs}  
Based on the expression of $\bar{\Sigma}_{d_0}$, the random variable $V_j$ or $Z_j =V_j/\sigma$ having the smallest value of fourth moments or kurtosis should be used. Under the additional condition $\esp \left[Z_1^2 \norme{\bo{Z} }^{2} \right] \geq 3 \rho_{\min}^2$, we can check that (see  Appendix \ref{app:theo:rtestfc})  
$$
\bar{\Sigma}_{d_0} \leq \frac{\left(\digamma_{d_0}^{\max} \right)^2}{N} \left[
 2\left(d\frac{\esp \left[Z_1^2 \norme{\bo{Z}}^{2} \right]}{3 \rho_{\min}^2}   \right)^{d_0} \indic_{d_0\leq d_0^*} + 2^d  
\left(\frac{\esp \left[Z_1^2 \norme{\bo{Z}}^{2} \right]}{3 \rho_{\min}^2}   \right)^{d_0} \indic_{d_0> d_0^*}
 \right] \, ,
$$
with $d_0^*:=\frac{(d-1)\ln(2)}{\ln(d)}$.    
\end{rem} 
 
In what follow, we are going to use $\esp\left[\left(\widehat{\M^c_{N, d_0}}(\bo{x}) -\M^c(\bo{x}) \right)^2\right] $ for the MSE of $\widehat{\M^c_{N, d_0}}$ and 
$$
\esp\left[\left(\widehat{\M^c_{N, d_0}} -\M^c \right)^2\right]  := \esp_{\bo{\X}}\left\{ \esp\left[\left(\widehat{\M^c_{N, d_0}}(\bo{\X}) -\M^c(\bo{\X}) \right)^2\right] \right\} \, , 
$$ 
for the expected or integrated MSE (IMSEs) of $\widehat{\M^c_{N, d_0}}$.  
    
\begin{corollary}  \label{coro:mse1} 
Given (\ref{eq:consttype1}), $r^*=d_0-1$, assume $\M \in \mathcal{H}_{\alpha}$ with $\alpha \in\left\{0, \max(d, d_0+2(L-d_0)) \right\}$; $\hh_k =\hh \to 0$ and (A1)-(A3) hold. Then, the MSE and IMSE share the same upper-bound given as follows: 
\begin{eqnarray}  
 \esp\left[\left(\widehat{\M^c_{N, d_0}} -\M^c \right)^2\right] & \leq &  2 D_{d_0, \rho_{\min}} \norme{\boldsymbol{\hh}^2}^{L-d_0}  \sigma^{2(L-d_0)}  K_{1,d_0-1, d_0}^{\max} \binom{d}{d_0} \left(\frac{1}{2  \rho_{\min}} \right)^{d_0}  +  D_{d_0, \rho_{\min}}^2  \nonumber \\
& &  + \bar{\Sigma}_{d_0} + 
 \norme{\boldsymbol{\hh}^2}^{2(L-d_0)}  \sigma^{4(L-d_0)}  \left(K_{1,d_0-1, d_0}^{\max}\right)^2  \binom{d}{d_0}^2 \left(\frac{1}{2  \rho_{\min}} \right)^{2 d_0}  \, .          \nonumber  
\end{eqnarray}        
 Moreover, if  $G_j=F_j$ with $j=1, \ldots, d$, then the IMSE is given by  
\begin{eqnarray}
\esp\left\{\esp\left[\left(\widehat{\M^c_{N, d_0}} -\M^c \right)^2\right] \right\} & \leq & \norme{\boldsymbol{\hh}^2}^{2(L-d_0)}  \sigma^{4(L-d_0)}  \left(K_{1,d_0-1, d_0}^{\max}\right)^2  \binom{d}{d_0}^2 \left(\frac{1}{2  \rho_{\min}} \right)^{2 d_0}     \nonumber \\
& & +  D_{d_0, \rho_{\min}}^2 + \bar{\Sigma}_{d_0}  \, .          \nonumber    
\end{eqnarray}          
\end{corollary}                      
\begin{preuve}
See Appendix \ref{app:coro:mse1}.  
\hfill $\square$          
\end{preuve}   

The presence of $D_{d_0, \rho_{\min}}$ in Corollary \ref{coro:mse1} is going to decrease the rates of convergence of our estimators without additional assumptions about $D_{d_0, \rho_{\min}}^2$. Corollary \ref{coro:mse11} starts giving such conditions and the associated MSEs and IMSEs.   

\begin{corollary}  \label{coro:mse11} 
Under the conditions of Corolary \ref{coro:mse1}, assume that $\M \in \mathcal{L}_{d_0, 0}$. Then, the upper-bound of the  IMSE and MSE is  given by   
$$  \norme{\boldsymbol{\hh}^2}^{2(L-d_0)}  \sigma^{4(L-d_0)}  \left(K_{1,d_0-1, d_0}^{\max}\right)^2  \binom{d}{d_0}^2 \left(\frac{1}{2  \rho_{\min}} \right)^{2 d_0}    
 +  \bar{\Sigma}_{d_0}   \, .  
$$  
Moreover, if  $V_k \sim \mathcal{U}(-\xi, \xi)$ with $\xi>0$ and $k=1, \ldots, d$, then this bound becomes 
$$
||\boldsymbol{\hh}^2||_1^{2(L-d_0)}  \xi^{4(L-d_0)}  \left( M_{d_0+ 2(L-d_0)} \Gamma_{d_0+2(L-d_0)}\right)^2 \binom{d}{d_0}^2 \left(\frac{1}{2  \rho_{\min}} \right)^{2d_0}  + \bar{\Sigma}_{d_0}  \, .   
$$  
\end{corollary}                     
\begin{preuve}
Using Corollary \ref{coro:biasfc1}, the results are straightforward.         
\hfill $\square$          
\end{preuve}   
 
Based on the upper-bounds of Corollary \ref{coro:mse11}, interesting choices of $\sigma$ or $\xi$ in one hand, and $\rho_{\min}$ and $\hh$ in the other hand help for obtaining the parametric rates of convergence due to the fact that $\bar{\Sigma}_{d_0}$ does not depend on $\hh$.  
 
\begin{corollary}  \label{coro:rate1} 
Let $r^*=d_0-1$, $L=d_0+1$. Assume $\M \in \mathcal{H}_{\alpha}$ with $\alpha \in \{0, \max(d, d_0+2)\}$; $\M \in \mathcal{L}_{d_0, 0}$  and (A1)-(A3) hold. If $\hh_k =\hh \propto N^{-\eta}$ with $\eta \in ]\frac{1}{4}, 1[$ and \\
 $\xi^2\leq \left(d  M_{d_0+ 2} \Gamma_{d_0+2} \binom{d}{d_0} \left(\frac{1}{2  \rho_{\min}} \right)^{d_0} \right)^{-1}$, then the upper-bound of MSE and IMSE is    
$$
\esp\left[\left(\widehat{\M^c_{N, d_0}} -\M^c \right)^2\right] \leq   
\frac{\left(\digamma_{d_0}^{\max} \right)^2}{N} \left[
 2\left(\frac{d(d+0.8)}{3 \rho_{\min}^2}   \right)^{d_0} \indic_{d_0\leq d_0^*} + 2^d  
\left(\frac{d+0.8}{3 \rho_{\min}^2}   \right)^{d_0} \indic_{d_0> d_0^*}
 \right] + \mathcal{O}(N^{-1})    \, ,         
$$
provided that $d+0.8 \geq 3 \rho_{\min}^2$. \\ 
Moreover, if $\rho_{\min} = \sqrt{\frac{c_0d(d+0.8)}{3}}$  with the real $1< c_0 \leq 2$, then    
$$
\esp\left[\left(\widehat{\M^c_{N, d_0}} -\M^c \right)^2\right] \leq  \frac{\left(\digamma_{d_0}^{\max} \right)^2}{N (c_0-1)} + \mathcal{O}(N^{-1}) \, .             
$$
\end{corollary}                    
\begin{preuve} 
See Appendix \ref{app:coro:rate1}.   
 \hfill $\square$              
\end{preuve} 
   
It is worth noting that the parametric rates of convergence are reached for any function that belongs to $\mathcal{L}_{d_0, 0}$ (see Corollary \ref{coro:rate1}). Moreover, taking $c_0 \in ]1,\, 2]$ leads to the dimension-free MSEs, which hold for particular distributions of the inputs given by  
$$
\X_j \sim F_j^*, \qquad  \rho_j^* \geq \rho_{\min}^* := \sqrt{\frac{c_0d(d+0.8)}{3}}, \; \forall\, j \in \{1, \ldots, d\} \,  .   
$$
For uniform distributions $\X_j \sim \mathcal{U}\left(a_j, b_j\right)$, we must have $b_j-a_j \leq  \frac{1}{\rho_{\min}^*}$. Obviously, such conditions are a bit strong, as a few distributions are covered. In the same sense, using $\widehat{\M^c_{N,d_0}}$ as an emulator of $\M^c$ for any sample point $\bo{x}$ of $\bo{\X} \sim F$ requires choosing $\X_j'$ such that its support contained that of $\X_j$, $j =1, \ldots, d$. Thus, the assumption (A3) given by $g_j>\rho_{\min}$ implicitly depends on the distribution of $\X_j$. For instance, given the bounded support of $\X_j'$, that is, $(a_j', \, b_j')$, we must have $\frac{1}{\rho_{\min}^*} >  b_j'-  a_j' \geq  b_j- a_j$ with $(a_j, \, b_j)$ the support of $\X_j$, limiting our ability to deploy $\widehat{\M^c_{N,d_0}}$ as a dimension-free, global emulator for some distributions of inputs. \\  
 
Nevertheless, the assumption $g_j \geq \rho_{\min}^*$ is always satisfied for the finite dimensionality $d$ if we are only interested in estimating $\M(\bo{\x}_0)$ for a given point $\bo{x}_0$, leading to local emulators. Indeed, taking $\X_j' \sim G_j$ to be depended on the point $\bo{x}^0$ at which $\M$ must be evaluated allows to enjoy the parametric rate of convergence and dimension-free MSEs for sample points falling in a neighborhood of $\bo{x}_0$. An example of such a choice is   
$$   
\X_j' \sim \mathcal{U}\left(x_{j}^0-\frac{1}{2\rho_{\min}^*}, x_{j}^0+ \frac{1}{2\rho_{\min}^*}\right) \, .
$$ 
However, different emulators are going to be built in order to estimate $\M(\bo{x})$ for any value $\bo{x}$ of $\bo{\X}$. Constructions of balls of given nodes and the radius $1/\rho_{\min}^*$ are an interesting perspective. 
  
\begin{rem} 
When $\M \notin \mathcal{L}_{d_0, 0}$, in-depth structural assumptions on $\M$ that should allow to enjoy the above MSEs
concern the truncation error, resulting from keeping only all the $|v|^{\mbox{h}}$ interactions or cross-partial derivatives with $|v| \leq d_0$. One way to handle it consists in choosing $d_0=d_{0, N}$ such that the residual bias is less than $1/N$ (i.e., $D_{d_{0, N}, \rho_{\min}}^2 < 1/N$) thanks to sensitivity indices.    
\end{rem}  
   
 While truncations are sometime necessary in higher-dimensions, it is interesting to have the rates of convergence without any truncation to cover lower or moderate dimensional functions for instance. 
         
\begin{corollary}  \label{coro:mse2} 
Let $d_0=d$; $r^* = d-1$, $L=d+1$ and $\hh_k=\hh$.   If $\M \in \mathcal{H}_{\alpha}$ with $\alpha \in \{0, d+1\}$, $\hh \to 0$, (A1)-(A3) hold; then the upper-bound of MSE and IMSE is   
\begin{equation}  
\esp\left[\left(\widehat{\M^c_{N}} -\M^c \right)^2\right] \leq    \sigma^2  \norme{\boldsymbol{\hh}}^2  M_{d+1}^2 K_{1:d}^2  \, \Gamma_{d+1}^2 \left(\frac{1}{2  \rho_{\min}} \right)^{2d} + \frac{\left(\digamma_{d}^{\max} \right)^2}{N}  \left[\left( \frac{\esp \left[Z_1^2 \norme{\bo{Z} }^{2} \right]}{3 \rho_{\min}^2} +1  \right)^{d}-1 \right]    \, .  \nonumber 
\end{equation}                  
 Moreover, if  $V_k \sim \mathcal{U}(-\xi, \xi)$ with $\xi>0$ and $k=1, \ldots, d$, then 
\begin{equation} \label{eq:upful2}
\esp\left[\left(\widehat{\M^c_{N}} -\M^c \right)^2\right] \leq    \xi^2 ||\boldsymbol{\hh}||_1^2   M_{d+1}^2 \, \Gamma_{d+1}^2 \left(\frac{1}{2  \rho_{\min}} \right)^{2d}
 +  \frac{\left(\digamma_{d}^{\max} \right)^2}{N}  \left[\left( \frac{d+0.8}{3 \rho_{\min}^2} +1  \right)^{d}-1 \right]   \, . \nonumber 
\end{equation}                 
\end{corollary} 
\begin{preuve}    
See Appendix \ref{app:coro:mse2}.  
\hfill $\square$         
\end{preuve}   
 
In the case of the full emulator of $\M^c$, remark that 
$$
\frac{\left(\digamma_{d}^{\max} \right)^2}{N}  \left[\left( \frac{d+0.8}{3 \rho_{\min}^2} +1  \right)^{d}-1 \right]
 \leq
\frac{d \left(\digamma_{d}^{\max} \right)^2(d+0.8)}{N\left(d(d+0.8) - 3 \rho_{\min}^2\right)}  \left[\left( d\frac{(d+0.8)}{3 \rho_{\min}^2}   \right)^{d}-1 \right] \, . 
$$
 
Based on Corollary \ref{coro:mse2}, different rates of convergence can be obtained depending on the the support of the input variables $\X_j \sim F_j$ via the choice of $\rho_{\min}$. 

\begin{corollary}  \label{coro:mse21} 
Let $r^*=d-1$ and  $L=d+1$.  Assume $\M \in \mathcal{H}_{\alpha}$ with $\alpha \in \{0,d+1\}$; 
 $\xi \leq \left[d   M_{d+1} \, \Gamma_{d+1} 
\left(\frac{1}{2  \rho_{\min}} \right)^{d} \right]^{-1}$; $\hh_k=\hh \propto N^{-\eta}$ with $\eta \in ]\frac{1}{2}, \, 1[$;  and (A1)-(A3) hold. Then, the (MSE and IMSE) rates of convergence are given as follows:  
$$
\esp\left[\left(\widehat{\M^c_{N}} -\M^c \right)^2\right] \leq 
\frac{\left(\digamma_{d}^{\max} \right)^2}{N}  \left[\left( \frac{d+0.8}{3 \rho_{\min}^2} +1  \right)^{d}-1 \right]
   +  \mathcal{O}\left(N^{-1}\right)\, .       
$$  
 
$\quad$ (i) If $\rho_{\min} = \sqrt{\frac{d+0.8}{3}}$, then  
$
\esp\left[\left(\widehat{\M^c_{N}} -\M^c \right)^2\right] \leq  \frac{\left(\digamma_{d}^{\max} \right)^2}{N}  (2^d-1)  + \mathcal{O}\left(N^{-1}\right)          
$.\\ 
  
$\quad$ (ii) If $\rho_{\min} = \sqrt{\frac{c_0d(d+0.8)}{3}}$ with $1 \leq c_0 \leq 2$, then     
\begin{equation}  
\esp\left[\left(\widehat{\M^c_{N}} -\M^c \right)^2\right] \leq  \frac{\left(\digamma_{d}^{\max} \right)^2}{N} \min\left\{ \frac{1}{c_0-1},\, \left( \frac{1}{c_0d} +1  \right)^{d}-1 \right\}  + \mathcal{O}\left(N^{-1}\right)  \nonumber  \, , 
\end{equation}     
with $\left( \frac{1}{c_0d} +1  \right)^{d}-1  \equiv \frac{1}{c_0}$ when $c_0d \to \infty$.
 \end{corollary}                
\begin{preuve} 
It is straightforward bearing in mind Corollaries \ref{coro:rate1}-\ref{coro:mse2}. 
     \hfill $\square$                                       
\end{preuve}         
    
Again, assumptions in Points (i)-(ii) are satisfied for fewer distributions, but they are  always satisfied if we are only interested in estimating $\M(\bo{\x}_0)$ for a given point $\bo{x}_0$, leading to build different local emulators. 

\subsubsection{Mean squared errors for every distribution of inputs} 
In this section, we are going to remove the assumption on the quasi-uniform distribution of inputs (A3) so as to cover any probability distribution of inputs. Note that (A3) is used to derive $\esp\left[E_k^2 \right] \leq \frac{1}{3 \rho_{\min}^2}$ and $\esp\left [ |E_k| \right] \leq \frac{1}{2 \rho_{\min}}$. Such an assumption can be avoided by using the following inequalities: 

$$
\esp\left [ |E_k| \right] =  \esp\left [\frac{\left[1- G(\X_k')\right]\indic_{X_k' \geq X_k} + G(\X_k') \indic_{X_k' < X_k}}{g_k(\X_k')} \right] \leq  \kappa_1 \, ,   
$$
with 
$$
\kappa_1 := \sup_{k \in \{1, \ldots, d\}} \sup_{x_k' \in \Omega_k}
\esp\left [\frac{F_k(\X_k')+ G(\X_k') - 2G(\X_k')F_k(\X_k')}{g_k(\X_k')} \right] \, ; 
$$
and   
$$
\esp\left [E_k^2 \right] =  \esp\left [ \frac{G^2(\X_k') + F_k(\X_k')-2G(\X_k')  F_k(\X_k')}{g_k^2(\X_k')} \right] \leq \kappa_2  \, ,   
$$
with
$$
\kappa_2 := \sup_{k \in \{1, \ldots, d\}} \sup_{x_k' \in \Omega_k} \esp\left [ \frac{G^2(x_k') + F_k(x_k')-2G(x_k')  F_k(x_k')}{g_k^2(x_k')} \right] \, .   
$$ 

Using such inequalities, the following results are straightforward keeping in mind Corollaries \ref{coro:rate1}-\ref{coro:mse21}. 
 
\begin{corollary}  \label{coro:rate1gen} 
Let $r^*=d_0-1$, $L=d_0+1$. Assume $\M \in \mathcal{H}_{\alpha}$ with $\alpha \in \{0, \max(d, d_0+2)\}$; $\M \in \mathcal{L}_{d_0, 0}$  and (A1)-(A2) hold. If $\hh_k =\hh \propto N^{-\eta}$ with $\eta \in ]\frac{1}{4}, 1[$ and \\
 $\xi^2\leq \left(d  M_{d_0+ 2} \Gamma_{d_0+2} \binom{d}{d_0} \kappa_1^{d_0} \right)^{-1}$, then we have 
$$
\esp\left[\left(\widehat{\M^c_{N, d_0}} -\M^c \right)^2\right] \leq   
\frac{\left(\digamma_{d_0}^{\max} \right)^2}{N} \left[
 2\left(\kappa_2 d(d+0.8) \right)^{d_0} \indic_{d_0\leq d_0^*} + 2^d  
\left(\kappa_2(d+0.8) \right)^{d_0} \indic_{d_0> d_0^*}
 \right]   \, , 
$$
provided that $\kappa_2(d+0.8) \geq 1$. 
\end{corollary} 

\begin{corollary}  \label{coro:mse21gen} 
Let $r^*=d-1$ and  $L=d+1$.  Assume $\M \in \mathcal{H}_{\alpha}$ with $\alpha \in \{0,d+1\}$; 
 $\xi \leq \left[d   M_{d+1} \, \Gamma_{d+1} 
\kappa_1^{d} \right]^{-1}$; $\hh_k=\hh \propto N^{-\eta}$ with $\eta \in ]\frac{1}{2}, \, 1[$;  and (A1)-(A2) hold. Then, 
$$
\esp\left[\left(\widehat{\M^c_{N}} -\M^c \right)^2\right] \leq 
\frac{\left(\digamma_{d}^{\max} \right)^2}{N}  \left[\left( \kappa_2 (d+0.8) +1  \right)^{d}-1 \right] \, . 
$$     
\end{corollary}           
Regarding the choice of $G_k$, it comes out from the expression 
$$
\esp\left [E_k^2 \right] =  \esp\left [ \frac{F_k(\X_k')+ G(\X_k') [G(\X_k') -2 F_k(\X_k')]}{g_k^2(\X_k')} \right] \, , 
$$
that  an interesting choice of $G_k$ must satisfy $G_k \leq  2 F_k$ in order to reduce the value of $\esp\left [E_k^2 \right]$. The following proposition gives interesting choices of $G_k$s. Recall that $\X_k \sim F_k$ and $\X_k' \sim G_k$ are supported on $\Omega_k,\, \Omega'_k$, respectively, with  $\Omega_k \subseteq \Omega'_k$.  
\begin{prop}  \label{prop:mixt}
Consider a PDF  $\rho'_k$ supported on $\Omega'_k \setminus \Omega_k$ and $\tau \in ]0,\, 1]$. The distribution $G_k$ defined as a mixture of $\rho_k$ and $\rho'_k$, that is, 
$$
g_k = \tau \rho_k\indic_{\Omega_k} +  (1-\tau) \rho_k' \indic_{\Omega_k' \setminus \Omega_k} \, , 
$$
allows to reduce $\esp\left [E_k^2 \right]$, provided that $\sup \{x: x \in \Omega\} \leq y_0$, $\forall\; y_0 \in \Omega_k' \setminus \Omega_k$.      
\end{prop}    
\begin{preuve}
We can check that $G_k \leq F_k$.      
   \hfill $\square$      
\end{preuve}
   
In what follows, we are going to consider $\tau=1$ (i.e., $G_k=F_k$) and $\tau < 1 $ with $\rho_k'$ the uniform distribution. 
								    	
\section{Illustrations} \label{sec:appli} 
In this section, we deploy our derivative-free emulators to approximate analytical functions. For the setting of different parameters needed, we rely on the results of Corollary \ref{coro:rate1}. Indeed, we use $V_k \sim \mathcal{U}(-\xi, \, \xi)$ with $\xi=\left(d \binom{d}{d_0} \left(\frac{1}{2  \rho_{\min}} \right)^{d_0} \right)^{-1/2}$; $L=d_0+1$ for the identified $d_0$. For each function, we use $NL$ runs of the model to construct our emulators, corresponding to $\M(\bo{\X}_i+\beta_\ell \hh\bo{V}_i)$ with 
\begin{itemize} 
\item $i=1, \ldots, N=500$;
\item $\beta_\ell \in \{0, \pm  2^{k-1}:  k= 1, \ldots, \frac{L-1}{2}\}$ if $L$ is odd and $\beta_\ell \in \{ \pm  2^{k}:  k=1, \ldots, \frac{L}{2}\}$ otherwise;  
\item $\hh = N^{-1}$.   
\end{itemize}
Then, we predict the output for $N=500$ sample points, that is, $\M(\bo{\X}_i)$. Finally,  sample values are generated using Sobol's sequence, and we compare the predictions associated with the initial distributions, that is, $G_k=F_k$ with a mixture distribution of $F_k$ (i.e., $\tau=0.9$) (see Proposition \ref{prop:mixt}).  
      
\subsection{Test~Functions} \label{sec:testfun} 
\subsubsection{Ishigami's Function ($d=3$)} 
\begin{figure}[!hbp]         
\begin{center}
\includegraphics[height=15cm,width=9cm,angle=270]{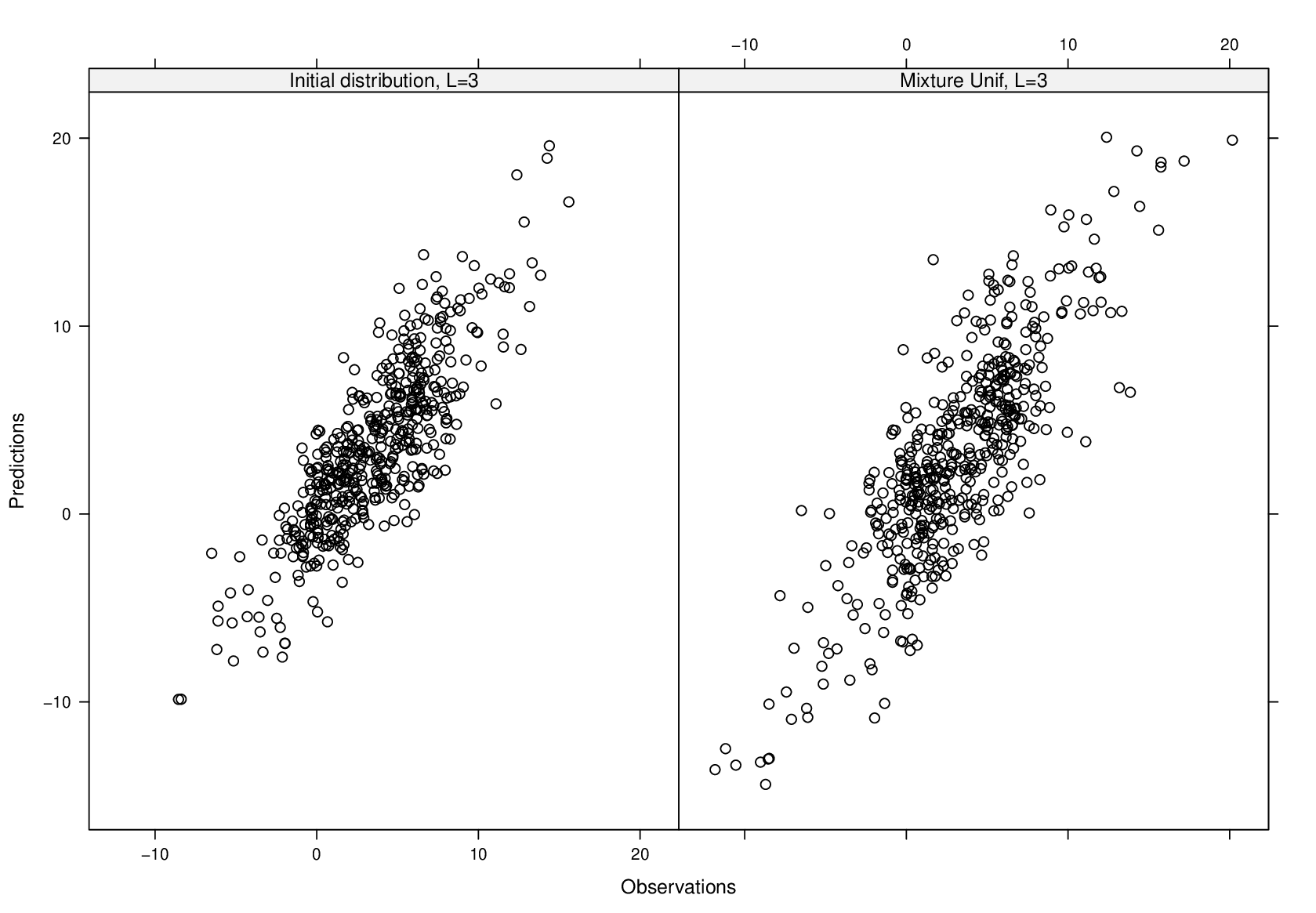}
\end{center}  
\caption{Predictions versus simulated outputs (observations) for Ishigami's function.}       
 \label{fig:ishi}              
\end{figure}

The Ishigami function is given by   
\begin{equation} \label{eq:ishi}
f(\bo{\x}) =  \sin(\x_1)+ 7\sin^2(\x_2) + 0.1 \x_3^4 \sin(\x_1) \, ,  \nonumber 
\end{equation}
with  $\X_j \sim \mathcal{U}(-\pi,\, \pi)$, $j=1, 2, 3$. The sensitivity indices are  $S_{1} = 0.3139$, $S_{2} = 0.4424$, $S_{3} = 0.0$,
$S_{T_1} = 0.567$, $S_{T_2} = 0.442$, and~$S_{T_3} = 0.243$. Thus, we have $d_0=2$ because $\sum_{j=1}^3 (S_j+ S_{T_j})=2.01$ (see Proposition \ref{prop:sedstrc}). Moreover, we are going to remove in our emulator the main effect term corresponding to $\X_3$, as $S_3=0$. Figure \ref{fig:ishi} depicts the predictions versus observations (i.e.,  simulated outputs) for $500$ sample points. We can see that our predictions are in line with the simulated values of the output for both distributions used.  
	                      
\subsubsection{Sobol's g-Function ($d=10$)}
The g-function (\cite{homma96}) is defined as follows:  
\begin{equation} \label{eq:gsobol}   
f(\bo{\x}) = \prod_{j=1}^{d=10}\frac{ |4 \x_j - 2| + a_j}{1 + a_j}\, .     \nonumber 
\end{equation}  
with $\X_j \sim \mathcal{U}(0,\, 1)$ $j=1, \ldots, d=10$. This function is differentiable almost everywhere and has different properties according to the values of $\bo{a} = (a_j, \, j = 1, \, 2, \, \ldots,\, d)$ (\cite{kucherenko09}). Indeed,    
\begin{itemize}
\item if $ \bo{a} = [0,\, 0,\, 6.52,\, 6.52,\, 6.52,\, 6.52,\, 6.52,\, 6.52,\, 6.52,\, 6.52]^\T$ (i.e., type A), we have $S_1 = S_2 = 0.39$, $S_j =0.0069$,   $\forall \, j > 2 $, $S_{T_1} = S_{T_2} = 0.54$,  and~$S_{T_j} = 0.013$,  $\forall \,  j > 2 $. We have $d_0=2$, as $\sum_{j=1}^d (S_j+S_{T_j})=2.01$. Moreover, we have $S_1+ S_{T_2} \approx 1$, suggesting to include only $\X_1$ and $\X_2$ in our emulator.  
\item  If $\bo{a} = [50,\, 50,\, 50,\, 50,\, 50,\, 50,\, 50,\, 50,\, 50,\, 50]^\T$ (i.e., type B), we have $S_j = S_{T_j} = 0.1$,  $\forall \, j \in  \{ 1,\, 2,\, \ldots,\, d\}$, leading to $d_0=1$.  
\item If $\bo{a} = [0,\, 0,\, 0,\, 0,\, 0,\, 0,\, 0,\, 0,\, 0,\, 0]^\T$ (i.e., type C),  we have $S_j = 0.02$ and  $ S_{T_j} = 0.27$,   $\forall\, j \in  \{ 1,\, 2,\, \ldots,\, d\}$, and we can check that $d_0=4$ and all the inputs are important. Thus, we have to include a lot of ANOVA components in our emulator with 
small effective effects since the variance of that function is fixed. More information is needed to better design the structure of this function.     
 \end{itemize}        
        
\begin{figure}[!hbp]         
\begin{center}
\includegraphics[height=15cm,width=9cm,angle=270]{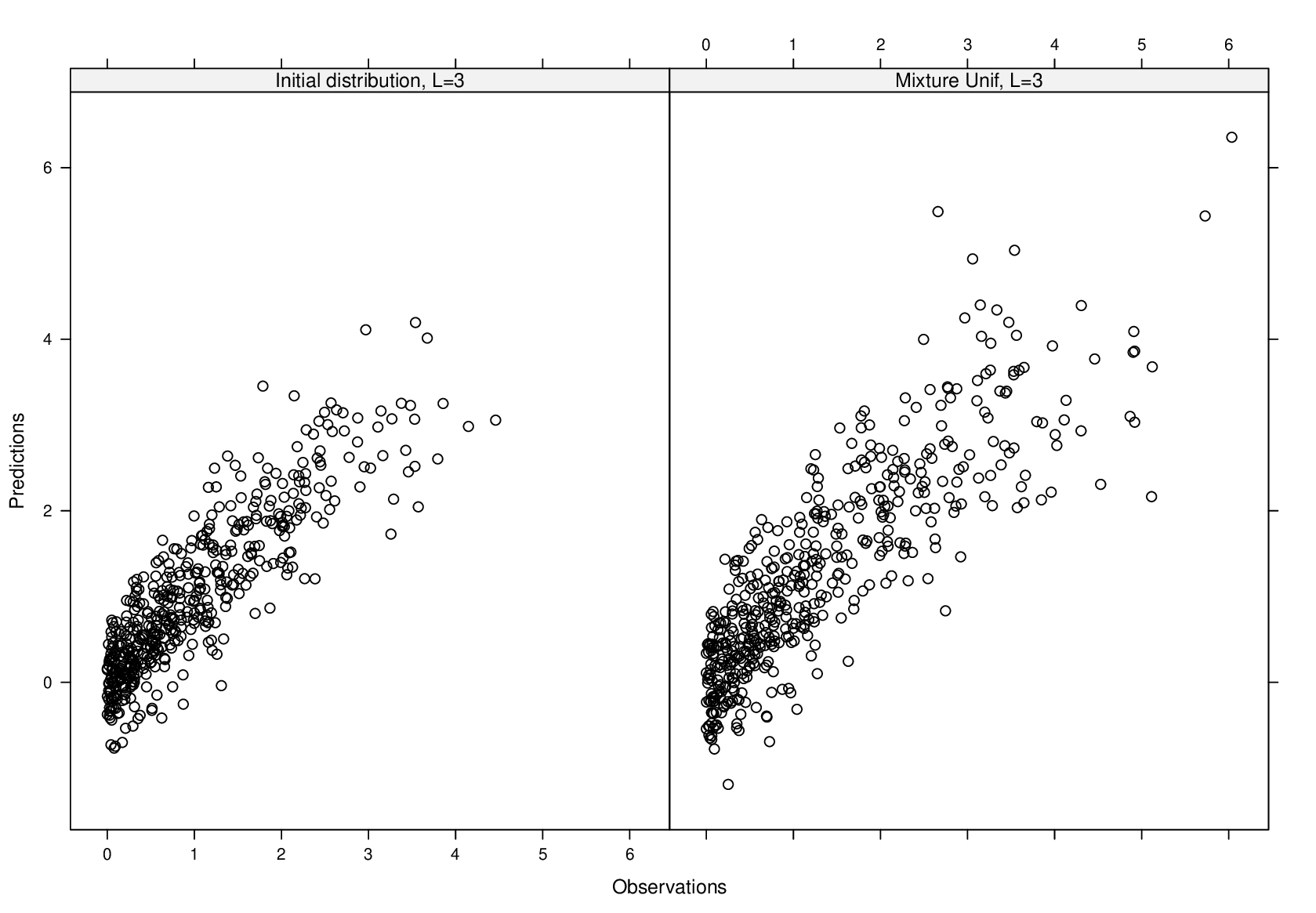}
\end{center} 
\caption{Predictions versus simulated outputs (observations) for the g-function of type A.}        
 \label{fig:sobla}          
\end{figure}

\begin{figure}[!hbp]         
\begin{center}
\includegraphics[height=15cm,width=9cm,angle=270]{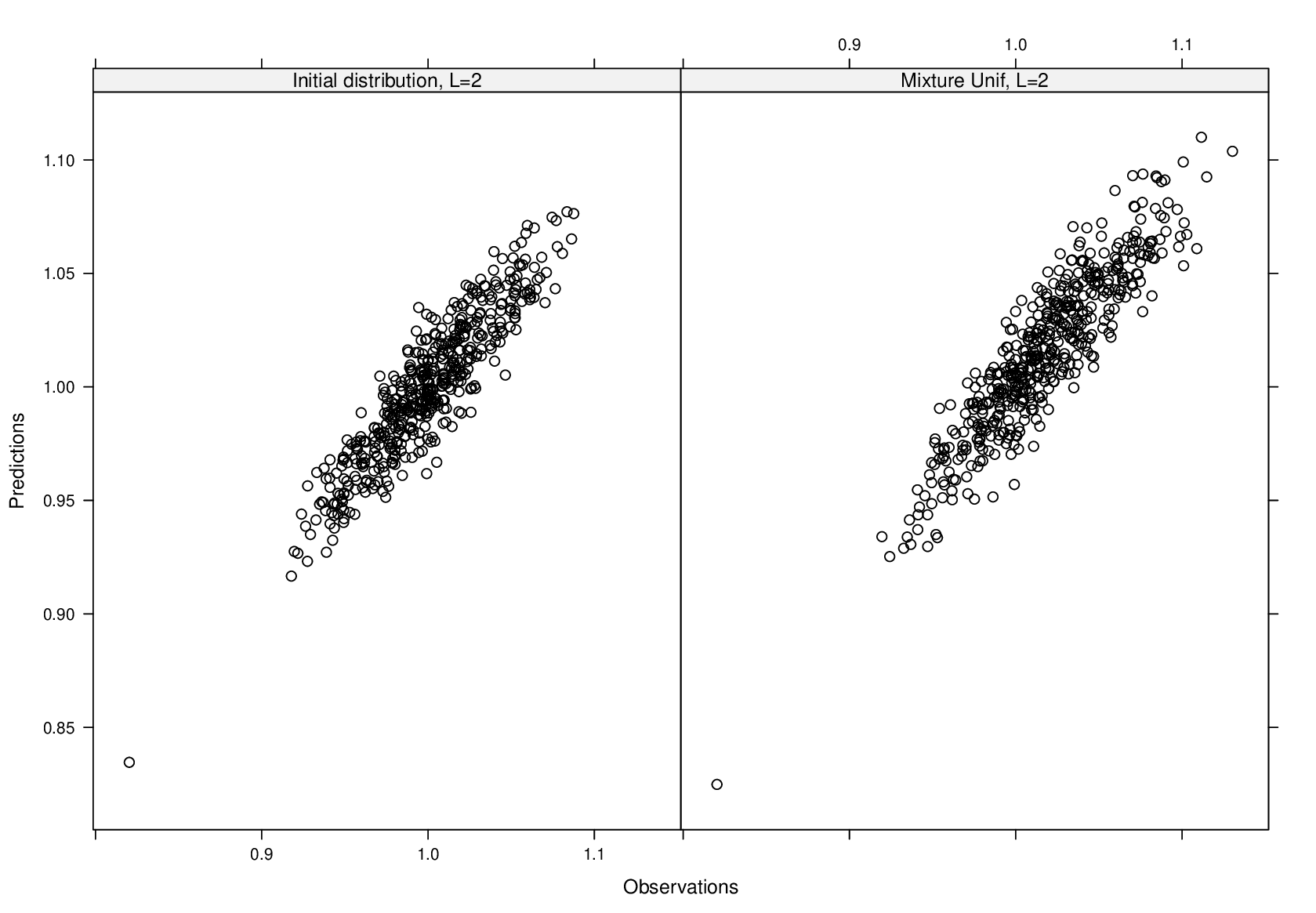}
\end{center} 
\caption{Predictions versus simulated outputs (observations) for the g-function of type B.}      
 \label{fig:sbolb}              
\end{figure}

Figures \ref{fig:sobla}-\ref{fig:sbolb} depict the predictions versus the simulated outputs (i.e., observations) of the g-function  of type A and type B, respectively, for $500$ sample points. While, we obtain quasi-perfect predictions in the case of the g-function of type B, those of type A face some difficulties in predicting small and null values.  

\section{Application: Heat diffusion models with stochastic initial conditions} \label{sec:realapp} 
We deploy our emulators to approximate a high-dimensional model defined by the one-dimensional ($1$-D) diffusion PDE with stochastic initial conditions, that is,  
$$  
\left\{ \begin{array}{lr}
\frac{\partial \M}{\partial t} - D \frac{\partial^2 \M }{\partial^2 x} =0,  & x \in ]0,\, 1[, \, t \in [0, T] \\
\M(x, t=0) =  Z(x)  & x \in [0, \, 1]\\
\M(x=0, t) =0, \quad  \M(x=1, t) =1,  &  t \in [0, \, T] \\    
\end{array} \right. \, ,      
$$   
where $D= 0.0011$ represents the diffusion coefficient. The quantity of interest (QoI) is given by $J(Z(x)) := \frac{1}{2} \int_0^T \int_0^{10} \left(\M(x, t)\right)^2 \, dxdt$. The spatial discretisation consists in subdividing the spatial domain $[0,\, 1]$ in $d$ equally-sized cells, leading to $d$ initial conditions or inputs, that is, $Z(x_j)$ with $j=1, \ldots, d$. We assume that the $d=50$ inputs are independent, and  $\X_j := Z(x_j) \sim \mathcal{U} \left(\sin(2\pi x_j) - 1.96,\;   \sin(2\pi x_j) +1.96 \right)$. A time step of $0.025$ is considered starting from $0$ up to $T=5$. \\  
It is knwon in \cite{lamboni24axioms} that the exact gradient can be computed as follows: 
$
\nabla_{Z}J(Z(x)) =  \M^a(x, 0)                     
$, 
where $\M^a(x, 0)$ stands for the adjoint model of $\M$ evaluated at $(x, t=0)$. Note that only one evaluation of such a function is needed to obtain the gradient of the QOI. The adjoint model is given by (see \cite{lamboni24axioms})
$$ 
\left\{ \begin{array}{lr}
-\frac{\partial \M^a}{\partial t} - D \frac{\partial^2 \M^a}{\partial^2 x} =\M,  & x \in ]0,\, 1[, \, t \in [0, T] \\
\M^a(x=0, t) =\M^a(x=1, t) = 0,    &  t \in [0, \, T]\\             
\M^a(x, T) = 0,   &  x \in [0, \, 1] \\           
\end{array} \right. \, .           
$$            
  
\begin{figure}[!hbp]               
\begin{center}
\includegraphics[height=15.25cm,width=16cm,angle=270]{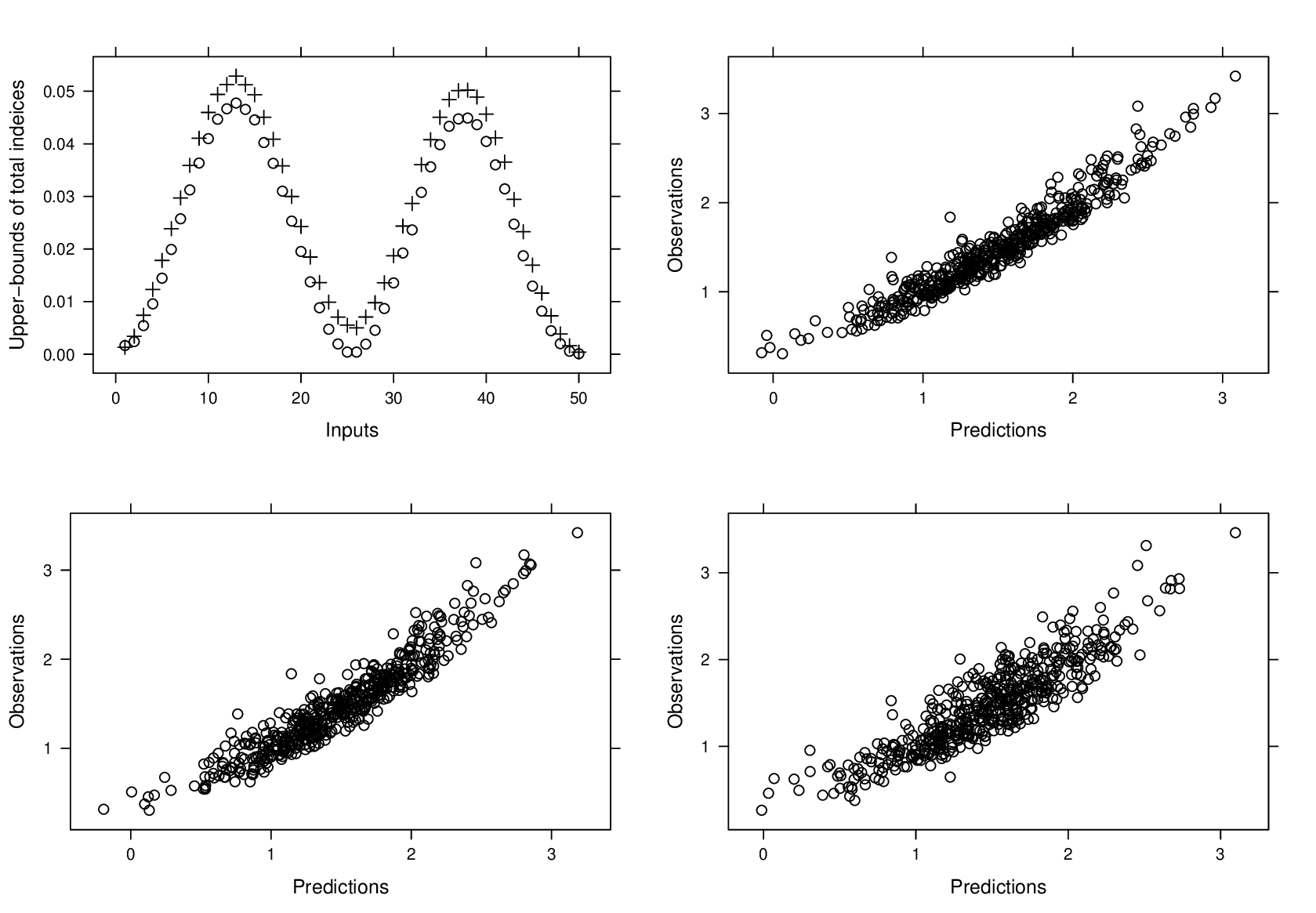}    
\end{center} 
\caption{Main indices ($\circ$) and upper-bounds of total indices ($+$) of $d=50$ inputs (left-top panel), and observations versus predictions using either derivative-based emulators (see right-top panel when including all component, and left-bottom panel otherwise) or derivative-free emulators (right-bottom panel).} 
 \label{fig:emulU} 
\end{figure}

The values of hyper-parameters derived in this paper and considered at the beginning of  Section \ref{sec:appli} are used to compute the results below. Using the exact values of the gradient (i.e., $\M^a(x, 0) $), we computed the main indices ($S_j$s) and the upper-bounds of the total indices ($U\!B_j$s) (see Figure \ref{fig:emulU}, left-top panel). It appears that the upper-bounds are almost equal to the main indices, showing the absence of interactions. This information is confirmed by the fact that 
 $\sum_{j=1}^{d=50} S_j =1.09$, leading to $d_0=1$. Based on this information, Figure \ref{fig:emulU} (right-top panel) depicts the predictions versus the observations (simulated outputs) using the derivative-based emulator with all the first-order partial derivatives. In the same sense, Figure \ref{fig:emulU} (left-bottom panel) depicts the observations versus predictions by including only the ANOVA components for which  $U\!B_j> 0.01$, that is, $37$ components. Both results are close together and are in line with the observations. Finally, Figure \ref{fig:emulU} (right-bottom panel) show the observations versus predictions for derivative-free emulators using only the components for which  $U\!B_j> 0.01$. It turns out that our emulators provide reliable estimations. As expected (see MSEs), the derivative-based emulator using exact values of derivatives performs better.  
           
\section{Conclusion} \label{sec:con}
In this paper, we have proposed simple, practical and sound emulators of simulators or estimators of functions using either available derivatives and distribution functions of inputs or derivative-free methods, such as stochastic approximations. Since our emulators  or estimators rely exactly on Db-ANOVA, Sobol' indices and their upper-bounds are  used to derive the appropriate structures of such emulators so as to reduce the epistemic uncertainty. The derivative-based and derivative-free emulators reach the parametric rate of convergence (i.e., $\mathsf{O}(N^{-1})$) and have dimension-free biases. Moreover, the former emulators enjoy the dimension-free MSEs when all cross-partial derivatives are available, and therefore, can cope with higher-dimensional models. But, the MSEs of the  derivative-free estimators depend on the dimensionality, and we have shown that the stability and accuracy of such emulators require about $N \approx (d+1)^d$ model runs  for full-emulators, and about $N \approx \min\left( d^{2d_0},\, 2^dd^{d_0} \right)$ runs for unbiased truncated emulators of order $d_0$.\\ 

To be able to deploy our emulators in practice, we have provided the best values of the hyper-parameters needed.
Numerical results have revealed that our emulators provide efficient predictions of models once the adequate structures of such models are used. While such results are promising, further improvements are going to be investigated in next works by i) considering distributions of $V_j$s that may help reducing the dimensionality in MSEs, ii) taking into account the discrepancy errors by using the output observations (rather than their mean only), and iii) considering local emulators. It is also worth investigating adaptations of such methods in the presence of empirical data.


\begin{appendices}     

\section{Derivations of unbiased truncations (Proposition \ref{prop:sedstrc})} \label{app:prop:sedstrc}
Keeping in mind the Sobol' indices, it is known that $ \sum_{\substack{v \subset \{1,\ldots, d\} \\ |v| >0}} S_v =1$, which comes down to 
$ \sum_{\substack{v \subset \{1,\ldots, d\} \\ |v|\leq d_0}} S_v =1$ for functions of the form
$
\M(\bo{\X}) = \sum_{\substack{v \subset \{1,\ldots, d\} \\  |v|\leq d_0}} \M_v(\bo{\X}_v )
$. Thus, we have 
$$
d_0 = \sum_{\substack{v \subset \{1,\ldots, d\} \\ |v|\leq d_0}} d_0 S_v \, ,
\qquad \qquad  
\sum_{j=1}^d S_{T_j} = \sum_{\substack{v \subset \{1,\ldots, d\} \\ |v|\leq d_0}} |v|  S_v
 \, .  
$$
Taking the difference yields  
$$
d_0 -\sum_{j=1}^d S_{T_j} = \sum_{\substack{v \subset \{1,\ldots, d\} \\ |v|\leq d_0 -1}} (d_0- |v|) S_v \Longrightarrow
d_0 -\sum_{j=1}^d  S_{T_j} - (d_0-1) \sum_{j=1}^d  S_{j}  = \sum_{\substack{v \subset \{1,\ldots, d\} \\ 2\leq |v|\leq d_0 -1}} (d_0- |v|) S_v \geq 0 
$$
which implies that $\sum_{j=1}^d  S_{T_j} + (d_0-1) \sum_{j=1}^d  S_{j} \leq d_0$. \\ 
Also, taking  
$
d_0- (d_0-1) \sum_{j=1}^d S_{T_j} =  \sum_{\substack{v \subset \{1,\ldots, d\} \\ |v|\leq d_0}} [d_0-(d_0-1)|v|]  S_v  
$ yields 
$$ 
d_0- (d_0-1) \sum_{j=1}^d S_{T_j} -  \sum_{j=1}^d  S_{j} = 
 \sum_{\substack{v \subset \{1,\ldots, d\} \\ 2\leq |v|\leq d_0}} [d_0-(d_0-1)|v|]  S_v  \leq 0 \, , 
$$
which implies $(d_0-1) \sum_{j=1}^d S_{T_j} +  \sum_{j=1}^d  S_{j} \geq d_0$.

\section{Proof of Theorem \ref{theo:sfall}} \label{app:theo:sfall}
Denote $\Vec{\boldsymbol{\imath}} :=(i_1, \ldots, i_d)  \in \N^d$, $||\Vec{\boldsymbol{\imath}}||_1 = i_1 +\ldots+ i_{d}$, $\Vec{\boldsymbol{\imath}}! := i_1! \ldots i_d!$, $\bo{\x}^{\Vec{\boldsymbol{\imath}}} := x_1^{i_1} \ldots x_d^{i_d}$, 
 and $\mathcal{D}^{(\Vec{\boldsymbol{\imath}})}\M := \partial^{i_1}_{1} \ldots \partial^{i_d}_{d} \M$. The Taylor expansion of  $\M\left(\bo{x} + \beta_\ell \boldsymbol{\hh}\bo{V} \right)$ about $\bo{x}$ of order $\alpha$ is given by    
\begin{eqnarray}                    
\M\left(\bo{x} + \beta_\ell \boldsymbol{\hh}\bo{V} \right)  &=&   \sum_{\rr = 0}^\alpha  \sum_{||\Vec{\boldsymbol{\imath}}||_1= \rr}  \frac{\mathcal{D}^{(\Vec{\boldsymbol{\imath}})}\M(\bo{x}) }{\Vec{\boldsymbol{\imath}}!} \beta_\ell ^\rr  \left( \boldsymbol{\hh}\bo{V} \right)^{\Vec{\boldsymbol{\imath}}} +\mathcal{O}\left(||\beta_\ell \boldsymbol{\hh}\bo{V} ||_1^{\alpha+1} \right) \, . \nonumber                   
\end{eqnarray}           
For any $w \subseteq \{1, \ldots, d\}$ with the cardinality $|w|$;  using $\indic_{w}(\cdot)$ for the indicator function  and $\Vec{\boldsymbol{w}} := \left(\indic_{w}(1), \ldots, \indic_{w}(d) \right)$ lead to $\mathcal{D}^{(\Vec{\boldsymbol{w}})}\M = \mathcal{D}^{|w|}\M$. Also, using $E_k :=  \frac{G_k \left(\X_k' \right) -\indic_{\X_k' \geq x_k }}{g_k\left(\X_k' \right)}$  implies that        
$
R_k = \frac{ V_k}{\hh_k \sigma^2} E_k              
$, $k=1, \ldots, d$.        \\  
Firstly, by evaluating the above expansion at $\bo{\X}'$ and taking the expectation w.r.t. $\bo{V}$, 
 $A := \sum_{\ell=1}^L C_{\ell}^{(|u|)} \esp_{\bo{V}} \left[ \M\left(\bo{\X}' + \beta_\ell \boldsymbol{\hh}\bo{V} \right) \bo{e}_{u}^{(1:d)}\left( \bo{R}\left(\bo{x}, \bo{\X}', \bo{V} \right) \right) \right]$ can be written as  
\begin{eqnarray}             
A &=&  \sum_{\rr\geq 0}  \sum_{||\Vec{\boldsymbol{\imath}}||_1= \rr}  \frac{\mathcal{D}^{(\Vec{\boldsymbol{\imath}})}\M(\bo{\X}')}{\Vec{\boldsymbol{\imath}}!}	\left( \sum_{\ell} C_{\ell}^{(|u|)}  \beta_{\ell}^{\rr} \right)  \esp_{\bo{V}} \left[\left( \boldsymbol{\hh}\bo{V} \right)^{\Vec{\boldsymbol{\imath}}}  
\sum_{\substack{w \subseteq \{1, \ldots, d\} \\ |w| =|u|}}	\left( \frac{ \bo{V}}{  \bo{\hh} \sigma^2}  \bo{E} \right)^{\Vec{\boldsymbol{w}}}	\right]    \nonumber  \\        
	&=& \sum_{\substack{\rr\geq 0 \\||\Vec{\boldsymbol{\imath}}||_1= \rr}}  
	\sum_{\substack{w \subseteq \{1, \ldots, d\} \\ |w| =|u|}}
	\frac{\mathcal{D}^{(\Vec{\boldsymbol{\imath}})}\M(\bo{\X}') \left( \sum_{\ell} C_{\ell}^{(|u|)}  \beta_{\ell}^{\rr} \right)}{\Vec{\boldsymbol{\imath}}! \, \sigma^{2|u|}} \esp_{\bo{V}} \left[\left(\bo{V} \right)^{\Vec{\boldsymbol{\imath}} +\Vec{\boldsymbol{w}}} \left(\boldsymbol{\hh} \right)^{\Vec{\boldsymbol{\imath}}-\Vec{\boldsymbol{w}}}	\left(\bo{E} \right)^{\Vec{\boldsymbol{w}}}	\right]  \,	\nonumber             
\end{eqnarray}         
We can see that $\esp_{\bo{V}} \left[ \left(\bo{V} \right)^{\Vec{\boldsymbol{\imath}} + \Vec{\boldsymbol{w}}}
\left( \boldsymbol{\hh} \right)^{\Vec{\boldsymbol{\imath}} -\Vec{\boldsymbol{w}}}	\right] \neq 0$ iff
$ \Vec{\boldsymbol{\imath}} + \Vec{\boldsymbol{w}} = 2\Vec{\bo{q}},\;  \forall \, \Vec{\bo{q}} \in \N^d$. 
Equation $\Vec{\boldsymbol{\imath}} + \Vec{\boldsymbol{w}} = 2\Vec{\bo{q}}$ implies  
$i_k =2q_k \geq 0$ if $k\notin w$ and $i_k=2q_k-1\geq 0$ otherwise. The last quantity is equivalent to 
 $i_k =2q_k +1$ when $k \in w$ with $q_k\in \N$, and it leads to $\Vec{\boldsymbol{\imath}} = 2\Vec{\bo{q}} + \Vec{\boldsymbol{w}}$, $\forall \, \Vec{\bo{q}} \in \N^d $, which also implies that $||  \Vec{\boldsymbol{\imath}}||_1 \geq ||\Vec{\boldsymbol{w}}||_1 = |u|$. \\ 
When $||\Vec{\bo{q}}||_1 =0$, we have $\Vec{\boldsymbol{\imath}} = \Vec{\boldsymbol{w}}$; $\mathcal{D}^{(\Vec{\boldsymbol{\imath}})}\M=\mathcal{D}^{w}\M$ and  
$
\esp \left[ \left(\bo{V} \right)^{\Vec{\boldsymbol{\imath}} + \Vec{\boldsymbol{w}}}
\left( \boldsymbol{\hh} \right)^{\Vec{\boldsymbol{\imath}} -\Vec{\boldsymbol{w}}}	\right]
= \esp \left[ \left(\bo{V} \right)^{2\Vec{\boldsymbol{w}}} \right] = \sigma^{2|u|}
$, by independence.  Thus, we have   

\begin{eqnarray}             
A &=&  \sum_{\substack{w \subseteq \{1, \ldots, d\} \\ |w| =|u|}} \mathcal{D}^{|w|}\M(\bo{\X}') \left( \sum_{\ell} C_{\ell}^{(|u|)}  \beta_{\ell}^{|u|} \right) \prod_{k \in w} E_k  + \sum_{\substack{s\geq 1 \\||\Vec{\bo{q}}||_1= s}}  \sum_{\substack{w \subseteq \{1, \ldots, d\} \\ |w| =|u|}}  \nonumber \\
 & & \times 
	\frac{\mathcal{D}^{(2\Vec{\bo{q}} + \Vec{\boldsymbol{w}})}\M(\bo{\X}') \left( \sum_{\ell} C_{\ell}^{(|u|)}  \beta_{\ell}^{|u|+2s} \right)}{(2\Vec{\bo{q}} + \Vec{\boldsymbol{w}})! \, \sigma^{2|u|}} \esp \left[\left(\bo{V} \right)^{2(\Vec{\boldsymbol{q}} +\Vec{\boldsymbol{w}})} \left(\boldsymbol{\hh} \right)^{2\Vec{\boldsymbol{q}}}	\left(\bo{E} \right)^{\Vec{\boldsymbol{w}}}	\right]  \,	, \nonumber            
\end{eqnarray}                      
using the change of variable $2s=\rr-|u|$. At this point, the setting $L=1, \, \beta_\ell=1$ and $C_{\ell}^{(|u|)}=1$ gives the approximation of   
$
\sum_{\substack{w \subseteq \{1, \ldots, d\} \\ |w| =|u|}} \mathcal{D}^{|w|}\M(\bo{\X}') \prod_{k \in w} \frac{G_k \left(\X_k' \right) -\indic_{\X_k' \geq x_k }}{g_k\left(\X_k' \right)}    
$ 
 of order $\mathcal{O}(\norme{\hh}^{2})$ for all $u \subseteq \{1, \ldots, d\}$. \\     
Secondly, taking the expectation w.r.t. $\bo{\X}'$ and the sum over $|u| =1, \ldots, d$, we obtain the result, that is,        
$$
\sum_{|u|=1}^d \esp_{\bo{X}'} \left[A \right] = \sum_{|u|=1}^d  
\sum_{\substack{w \subseteq \{1, \ldots, d\} \\ |w| =|u|}} \esp_{\bo{\X}'}\left[ \mathcal{D}^{|w|}\M(\bo{\X}') \prod_{k \in w} \frac{G_k \left(\X_k' \right) -\indic_{\X_k' \geq x_k }}{g_k\left(\X_k' \right)} \right]   + \mathcal{O}\left(\norme{\hh}^{2}\right) \, ,           
$$    
bearing in mind Equation (\ref{eq:deriGene}).  \\ 
Thirdly,  we can increase this order up to $\mathcal{O}\left( \norme{\bo{\hh}}^{2L} \right)$ by using the constraints $\sum_{\ell=1}^L C_{\ell}^{(|u|)}  \beta_{\ell}^{2s+|u|} =\delta_{0, s}$ $s=0, 1, \ldots, L$ to eliminate some higher-order terms. Thus, the order $\mathcal{O}(\norme{\hh}^{2L})$ is reached. Other constraints can lead to the orders $\mathcal{O}(\norme{\hh}^{2 \alpha_{|u|}})$ with $\alpha_{|u|} =1, \ldots, L$. 

\section{Proof of Lemma \ref{lem:trucbias}} \label{app:lem:trucbias}
Recall that  $\left| E_k \right| =  \left|\frac{G_k \left(\X_k' \right) -\indic_{\X_k' \geq x_k }}{g_k\left(\X_k' \right)} \right|$. By using the definition of the absolute value and the fact that $0\leq G(x)  \leq 1$, we can check that 
$\esp \left| E_k \right|  \leq \frac{1}{2  \rho_{\min}}$. Also, using the following inequality (see Lemma 1 in \cite{bach16}) 
$$ 
M_{|w|}' \leq  2 \gamma_0   \left(M_{0}'\right)^{1-|w|/d}   \left( M_{d}'\right)^{|w|/d} 
= 2 \gamma_0 M_{0}' \left(\frac{M_{d}'}{M_{0}'} \right)^{|w|/d} \, ,    
$$ 
for a given $\gamma_0$, we can write   
\begin{eqnarray}     
D_{d_0, \rho_{\min}} & :=& \sum_{\substack{w \subseteq \{1, \ldots, d\} \\ |w| >d_0}} M_{|w|}'  \left(\frac{1}{2  \rho_{\min}} \right)^{|w|} \leq  2 \gamma_0 M_{0}' \sum_{\substack{w \subseteq \{1, \ldots, d\} \\ |w| >d_0}}  \left[\frac{1}{2  \rho_{\min}} \left(\frac{M_{d}'}{M_{0}'} \right)^{1/d} \right]^{|w|} \nonumber \\
&=&   2 \gamma_0 M_{0}'\left\{ \sum_{\substack{w \subseteq \{1, \ldots, d\}}}  \left[\frac{1}{2  \rho_{\min}} \left(\frac{M_{d}'}{M_{0}'} \right)^{1/d} \right]^{|w|} - \sum_{\substack{w \subseteq \{1, \ldots, d\} \\ |w| \leq d_0}}  \left[\frac{1}{2  \rho_{\min}} \left(\frac{M_{d}'}{M_{0}'} \right)^{1/d} \right]^{|w|} \right\} \nonumber \\
& \leq & 2 \gamma_0 M_{0}'\left\{  
\sum_{\substack{\ell=0}}^d  \binom{d}{\ell} \left[\frac{1}{2  \rho_{\min}} \left(\frac{M_{d}'}{M_{0}'} \right)^{1/d} \right]^{\ell} - \sum_{\substack{\ell=0}}^{d_0}  \binom{d_0}{\ell} \left[\frac{1}{2  \rho_{\min}} \left(\frac{M_{d}'}{M_{0}'} \right)^{1/d} \right]^{\ell} \right\} \nonumber \\     
&= & 2 \gamma_0 M_{0}'\left\{  
\left[\frac{1}{2  \rho_{\min}} \left(\frac{M_{d}'}{M_{0}'} \right)^{1/d} +1\right]^{d} -
\left[\frac{1}{2  \rho_{\min}} \left(\frac{M_{d}'}{M_{0}'} \right)^{1/d} +1\right]^{d_0}
\right\} \nonumber \\        
&=& 2 \gamma_0 M_{0}' \left[\frac{1}{2  \rho_{\min}} \left(\frac{M_{d}'}{M_{0}'} \right)^{1/d} +1 \right]^{d_0}
\left\{ \left[\frac{1}{2  \rho_{\min}} \left(\frac{M_{d}'}{M_{0}'} \right)^{1/d} +1 \right]^{d-d_0} -1 \right\} 
  \nonumber    \, ,   
\end{eqnarray}   
and the result holds.  
   
\section{Proof of Theorem \ref{theo:biasfc}} \label{app:theo:biasfc}
Recall that $R_k = \frac{ V_k}{\hh_k \sigma^2} E_k$ with $E_k :=  \frac{G_k \left(\X_k' \right) -\indic_{\X_k' \geq x_k }}{g_k\left(\X_k' \right)}$, $k=1, \ldots, d$.  Using  
$$
A_1  :=   \sum_{p=1}^{d_0} \sum_{\ell=1}^{L} C_{\ell}^{(p)} \esp\left\{ \left[\M\left(\bo{\X}' + \beta_\ell \boldsymbol{\hh}\bo{V} \right) \right] \bo{e}_{p}^{(1:d)}\left( \bo{R}\left(\bo{x}, \bo{\X}', \bo{V} \right) \right) \right\} \, ,   
$$ 
 the bias $B := A_1-\M^c(\bo{x})$ becomes       
\begin{eqnarray}
B &=& \sum_{\substack{w \subseteq \{1, \ldots, d\} \\ 0< |w| \leq  d_0}}  \left\{  \sum_{\ell=1}^{L} C_{\ell}^{(|w|)} \esp\left[ \M\left(\bo{\X}' + \beta_\ell \boldsymbol{\hh}\bo{V} \right)  \prod_{k \in w} \frac{ V_k}{\hh_k \sigma^2} E_k   \right] - \esp_{\bo{\X}'}\left[ \mathcal{D}^{|w|}\M(\bo{\X}') \prod_{k \in w} E_k \right] \right\}  \nonumber\\ 
& & - \sum_{\substack{w \subseteq \{1, \ldots, d\} \\ |w| >d_0}} \esp_{\bo{\X}'}\left[ \mathcal{D}^{|w|}\M(\bo{\X}') \prod_{k \in w} E_k \right]  \nonumber\\
&=&   \sum_{\substack{w \subseteq \{1, \ldots, d_0\} \\ |w| >0}} \esp_{\bo{\X}'} \left\{ \left(\prod_{k \in w} E_k \right) \left( \sum_{\ell=1}^{L} C_{\ell}^{(|w|)} \esp_{V}\left[ \M\left(\bo{\X}' + \beta_\ell \boldsymbol{\hh}\bo{V} \right)  \prod_{k \in w} \frac{ V_k}{\hh_k \sigma^2}  \right] -  \mathcal{D}^{|w|}\M(\bo{\X}')  \right) \right\} \nonumber\\ 
& & - \sum_{\substack{w \subseteq \{1, \ldots, d\} \\ |w| >d_0}} \esp_{\bo{\X}'}\left[ \mathcal{D}^{|w|}\M(\bo{\X}') \prod_{k \in w} E_k \right]  \nonumber \, .    
\end{eqnarray}    
Note that the quantity $\sum_{\ell=1}^{L} C_{\ell}^{(|w|)} \esp_{V}\left[ \M\left(\bo{\X}' + \beta_\ell \boldsymbol{\hh}\bo{V} \right)  \prod_{k \in w} \frac{ V_k}{\hh_k \sigma^2}  \right] -  \mathcal{D}^{|w|}\M(\bo{\X}') $ has been investigated in \cite{lamboni24stats}. To make use of such results in our context given by Equation (\ref{eq:consttype1}), let  
$   
L'_w := \left(\left[\frac{r^*-|w|}{2} \right] + L-r^* \right) \indic_{|w|\leq r^*} +  (L-r^*-1) \indic_{|w|> r^*}
$. Thus, we have  
\begin{equation}  \label{eq:bias0}                          
\displaystyle    
  \left|
\mathcal{D}^{|w|}\M(\bo{x}) - \sum_{\ell=1}^{L} C_{\ell}^{(|w|)} \, \esp \left[ \M\left(\bo{x} + \beta_\ell \boldsymbol{\hh}\bo{V} \right) \prod_{k \in w} \frac{ V_k}{\hh_k \sigma^2} \right]
 \right| \leq \sigma^{2L'_w }  M_{|w|+2L'_w } K_{1,L'_w}  \norme{\boldsymbol{\hh}^2}^{L'_w} \, .         
\end{equation}           
 When $V_k \sim \mathcal{U}(-\xi, \xi)$ with $\xi>0$ and $k=1, \ldots,d$, then
\begin{equation}  \label{eq:bias20}                  
\displaystyle            
   \left| 
\mathcal{D}^{|w|}\M(\bo{x}) - \sum_{\ell=1}^{L} C_{\ell}^{(|w|)} \, \esp \left[ \M\left(\bo{x} + \beta_\ell \boldsymbol{\hh}\bo{V} \right) \prod_{k \in w} \frac{ V_k}{\hh_k \sigma^2} \right] 
 \right| \leq   M_{|w|+2L'_w } \xi^{2L'_w }  ||\boldsymbol{\hh}^2||_1^{L'_w}   \Gamma_{|w|+2L'_w}   \, . 
\end{equation}   
 Using Equation (\ref{eq:bias0}); $g_k > \rho_{\min}$  and the fact that 
 $\esp \left| E_k \right| \leq \frac{1}{2  \rho_{\min}}$, we can write          
\begin{eqnarray}   
|B|  &\leq& \sum_{\substack{w \subseteq \{1, \ldots, d\} \\ 0< |w| \leq  d_0}}   \sigma^{2L'_w  }  M_{|w|+2L'_w } K_{1,L'_w }  \norme{\boldsymbol{\hh}^2}^{L'_w } \prod_{k \in w} \esp\left[\left| E_k \right| \right]  + \sum_{\substack{w \subseteq \{1, \ldots, d\} \\ |w| >d_0}} M_{|w|}' \prod_{k \in w} \esp\left[\left| E_k \right| \right]  \nonumber \\    
& \leq &  \sum_{\substack{w \subseteq \{1, \ldots, d\} \\ 0< |w| \leq  d_0}}   \sigma^{2L'_w  }  M_{|w|+2L'_w } K_{1,L'_w }  \norme{\boldsymbol{\hh}^2}^{L'_w } \prod_{k \in w} \frac{1}{2  \rho_{\min}}  + \sum_{\substack{w \subseteq \{1, \ldots, d\} \\ |w| >d_0}} M_{|w|}'  \prod_{k \in w} \frac{1}{2  \rho_{\min}}   \nonumber  \, , 
\end{eqnarray}          
where $M_{|w|}' = \norminf{\mathcal{D}^{|w|}\M}$. The results hold by using Lemma \ref{lem:trucbias}.

\section{Proof of Corollary \ref{coro:biasfc1}} \label{app:coro:biasfc1}
Firstly, keeping in mind Theorem \ref{theo:biasfc}, we can see that the smallest value of $L'_w$ is $L-r^*-1$, which is reached when $|w|>r^*$. Thus, the bias verifies  
\begin{eqnarray}  
\left| B \right| &\leq &   \norme{\boldsymbol{\hh}^2}^{L-r^*-1}
  \sigma^{2(L-r^*-1)}  
\sum_{\substack{w \subseteq \{1, \ldots, d\} \\ r^* < |w| \leq  d_0}} K_{w, (L-r^*-1)}  M_{|w|+ 2(L-r^*-1) }  \left(\frac{1}{2  \rho_{\min}} \right)^{|w|} \nonumber \\ 
& & +  D_{d_0, \rho_{\min}} +   \mathcal{O}\left(\norme{\boldsymbol{\hh}^2}^{L-r^*-1} \right) \, .  \nonumber 
\end{eqnarray}                     
Secondly, using $K_{1,r^*, d_0}^{\max}$, we can write   
\begin{eqnarray}
A_3 &:=& \sum_{\substack{w \subseteq \{1, \ldots, d\} \\ r^* < |w| \leq  d_0}} K_{w, (L-r^*-1)}  M_{|w|+ 2(L-r^*-1) }  \left(\frac{1}{2  \rho_{\min}} \right)^{|w|}  \leq  K_{1,r^*, d_0}^{\max} \sum_{\imath =r^*+1}^{d_0} \binom{d}{\imath} \left(\frac{1}{2  \rho_{\min}} \right)^{\imath} \nonumber \\  
&=&  2 K_{1, d_0}^{\max}    \rho_{\min}
\left(\frac{d}{2  \rho_{\min}}\right)^{r^*+1} \frac{ \left(\frac{d}{2  \rho_{\min}}\right)^{d_0-r^*} - 1}{d-2  \rho_{\min}} \, ,  \nonumber    
\end{eqnarray}          
because $ \binom{d}{\imath} \leq d ^{\imath}$ and   
$
\sum_{\imath =r^*+1}^{d_0} \binom{d}{\imath} \left(\frac{1}{2  \rho_{\min}} \right)^{\imath} \leq   
\sum_{\imath =r^*+1}^{d_0} \left(\frac{d}{2  \rho_{\min}}\right)^{\imath} =  \left(\frac{d}{2  \rho_{\min}}\right)^{r^*+1} \frac{ \left(\frac{d}{2  \rho_{\min}}\right)^{d_0-r^*} - 1}{\frac{d}{2  \rho_{\min}}  -1}
$.    \\
Finally, if  $V_k \sim \mathcal{U}(-\xi, \xi)$ with $\xi>0$ and $k=1, \ldots, d$, the following bias is used to derve the result:    
\begin{eqnarray}
\left| B \right| &\leq &   ||\boldsymbol{\hh}^2||_1^{L-r^*-1}  \xi^{2(L-r^*-1)}  \sum_{\substack{w \subseteq \{1, \ldots, d\} \\ r^*< |w| \leq  d_0}}   
 M_{|w|+ 2(L-r^*-1)} \Gamma_{|w|+2(L-r^*-1)}  \left(\frac{1}{2  \rho_{\min}} \right)^{|w|}  \nonumber  \\ 
& &  + D_{d_0, \rho_{\min}} + \mathcal{O}\left( ||\boldsymbol{\hh}^2||_1^{L-r^*-1} \right)   \, .  \nonumber    
\end{eqnarray}       

\section{Proof of Corollary \ref{coro:biasfc2}} \label{app:coro:biasfc2}
For $r^*< |w|\leq d-1$, we can see that $L'_w =1$, and the order of approximation in Corollary \ref{coro:biasfc1} becomes $\mathcal{O}\left(\norme{\boldsymbol{\hh}^2} \right)$  or $\mathcal{O}\left( ||\boldsymbol{\hh}^2||_1 \right)$ because  $|w|+2\leq d+1$ and $L=r^*+2$. When $|w|=d$, the smallest order is obtained thanks to [\cite{lamboni24stats}, Corollary 2].   

\section{Proof of Theorem \ref{theo:rtestfc}} \label{app:theo:rtestfc}
For the variance of our emulator,  we can write      
\begin{eqnarray}    
\var\left[\widehat{\M^c_{N, d_0}}(\bo{x}) \right]  &= & \frac{1}{N} \var \left[ \sum_{p=1}^{d_0} \sum_{\ell=1}^{L} C_{\ell}^{(p)} \M\left(\bo{\X}' + \beta_\ell \boldsymbol{\hh}\bo{V} \right) \bo{e}_{p}^{(1:d)}\left( \bo{R}\left(\bo{x}, \bo{\X}', \bo{V} \right) \right) \right] \nonumber \\  
&\leq & 
\frac{1}{N} \esp \left[ \left\{\sum_{p=1}^{d_0} \sum_{\ell=1}^{L} C_{\ell}^{(p)} \M\left(\bo{\X}' + \beta_\ell \boldsymbol{\hh}\bo{V} \right) \bo{e}_{p}^{(1:d)}\left( \bo{R}\left(\bo{x}, \bo{\X}', \bo{V} \right) \right) \right\}^2 \right] \nonumber \\       
& \leq &    \frac{1}{N}  \esp \left[\left\{\sum_{p=1}^{d_0} 
\sum_{\substack{w \subseteq \{1, \ldots, d\} \\  |w|=p}} \left(\prod_{k \in w} \frac{ V_k}{\hh_k \sigma^2} E_k \right)  \left(
 \sum_{\ell=1}^{L} C_{\ell}^{(p=|w|)} \M\left(\bo{\X}' + \beta_\ell \boldsymbol{\hh}\bo{V} \right) \right) \right\}^2 \right]  \nonumber  \, .           
\end{eqnarray} 
Using $A_2 :=  \sum_{\ell=1}^{L} C_{\ell}^{(|w|)} \M\left(\bo{\X}' + \beta_\ell \boldsymbol{\hh}\bo{V} \right)$, and keeping in mind Equation (\ref{eq:consttype1}), we can write 
\begin{eqnarray} 
\left| A_2 \right| &=&  \left|\sum_{\ell=1}^{L} C_{\ell}^{(|w|)} \left[ \M\left(\bo{\X}' + \beta_\ell \boldsymbol{\hh}\bo{V} \right) - 
\sum_{\rr = 0}^{\min(r^*, |w|-1)}  \sum_{||\Vec{\boldsymbol{\imath}}||_1= \rr}  \frac{\mathcal{D}^{(\Vec{\boldsymbol{\imath}})}\M(\bo{\X}') }{\Vec{\boldsymbol{\imath}}!} \beta_\ell ^\rr  \left( \boldsymbol{\hh}\bo{V} \right)^{\Vec{\boldsymbol{\imath}}} \right] \right| \nonumber \\
& \leq &  M_{\min(r^*, |w|-1) +1} \norme{\boldsymbol{\hh}\bo{V} }^{\min(r^*, |w|-1) +1}  \sum_{\ell=1}^{L} \left| C_{\ell}^{(|w|)} \beta_\ell^{\min(r^*, |w|-1) +1}  \right| \, ,   \nonumber 
\end{eqnarray}   
because $\M \in \mathcal{H}_{|w|+1}$.  Keeping in mind that $\Gamma_{\min(r^*, |w|-1) +1} = \sum_{\ell=1}^{L} \left| C_{\ell}^{(|w|)} \beta_\ell^{\min(r^*, |w|-1) +1}  \right|$ and using $\digamma_{r^*,|w|} :=   M_{\min(r^*, |w|-1) +1} \Gamma_{\min(r^*, |w|-1) +1}$, the variance becomes  
 
\begin{eqnarray}    
\var\left[\widehat{\M^c_{N, d_0}}(\bo{x}) \right]  
& \leq &    \frac{1}{N} \esp \left[\left\{
\sum_{\substack{w \subseteq \{1, \ldots, d\} \\  0<|w| \leq d_0}}  \digamma_{r^*,|w|} \norme{\boldsymbol{\hh}\bo{V} }^{\min(r^*, |w|-1) +1} \left(\prod_{k \in w} \frac{ V_k}{\hh_k \sigma^2} E_k \right)  \right\}^2 \right]  \nonumber  \\ 
 & \leq &   \frac{\left(\digamma_{r^*,d_0}^{\max} \right)^2}{N} \esp \left[\left\{
\sum_{\substack{w \subseteq \{1, \ldots, d\} \\ 0< |w| \leq d_0}} \prod_{k \in w} \left( \frac{ V_k  E_k}{\hh_k \sigma^2}  \norme{\boldsymbol{\hh}\bo{V} }^{\frac{\min(r^*, |w|-1) +1}{|w|}} \right)\right\}^2 \right]  \nonumber  \\    
 & = &    \frac{\left(\digamma_{r^*,d_0}^{\max} \right)^2}{N} \esp \left[ 
\sum_{\substack{w \subseteq \{1, \ldots, d\} \\ 0< |w| \leq d_0}}
 \prod_{k \in w} \left( \frac{ V_k  E_k}{\hh_k \sigma^2}  \norme{\boldsymbol{\hh}\bo{V} }^{\frac{\min(r^*, |w|-1) +1}{|w|}} \right)^2  \right]  \nonumber \, ,  
\end{eqnarray}   
because $\esp[V_k] =0$,  $\esp\left[V_k\norme{\boldsymbol{\hh}\bo{V} }^{\frac{\min(r^*, |w|-1) +1}{|w|}}\right] =0$ and when $w \neq w'$,         
$$
\esp\left[\prod_{k \in w} \left( \frac{ V_k  E_k}{\hh_k \sigma^2}  \norme{\boldsymbol{\hh}\bo{V} }^{\frac{\min(r^*, |w|-1) +1}{|w|}} \right)  
 \prod_{\ell \in w'} \left( \frac{ V_\ell  E_\ell}{\hh_\ell \sigma^2}  \norme{\boldsymbol{\hh}\bo{V} }^{\frac{\min(r^*, |w|-1) +1}{|w|}} \right)  \right] =0 \,. 
$$ 
For the second result, by expanding $E_k^2$ and knowing that $0\leq G(x) \leq 1$, we can see that  $\esp[E_k^2] \leq \frac{1}{3\rho_{\min}^2}$.  Also, using $Z_k=V_k/\sigma$ and the fact that $\hh_k=\hh,\, r^*=d_0-1$ and $|w|\leq d_0$, we have      
\begin{eqnarray}    \label{eq:intub}
\var\left[\widehat{\M^c_{N, d_0}}(\bo{x}) \right]    
 & \leq &    \frac{\left(\digamma_{d_0}^{\max} \right)^2}{N} 
\sum_{\substack{w \subseteq \{1, \ldots, d\} \\ 0< |w| \leq d_0}}
 \prod_{k \in w} \left( \frac{\esp \left[V_1^2 \norme{\bo{V} }^{2} \right]}{\sigma^4} \esp[E_k^2]\right)  \nonumber   \\    
& \leq &    \frac{\left(\digamma_{d_0}^{\max} \right)^2}{N} 
\sum_{\substack{w \subseteq \{1, \ldots, d\} \\ 0< |w| \leq d_0}}
 \prod_{k \in w} \left(\esp \left[Z_1^2 \norme{\bo{Z} }^{2} \right] \esp[E_k^2] \right)    \\      
 & \leq &    \frac{\left(\digamma_{d_0}^{\max} \right)^2}{N} 
\sum_{\substack{w \subseteq \{1, \ldots, d\} \\ 0< |w| \leq d_0}}
\left( \frac{\esp \left[Z_1^2 \norme{\bo{Z} }^{2} \right]}{3 \rho_{\min}^2}   \right)^{|w|}
=  \frac{\left(\digamma_{d_0}^{\max} \right)^2}{N} 
\sum_{p=1}^{d_0} \binom{d}{p} \left( \frac{\esp \left[Z_1^2 \norme{\bo{Z} }^{2} \right]}{3 \rho_{\min}^2}   \right)^{p}   \nonumber \\ 
&\leq & \frac{\left(\digamma_{d_0}^{\max} \right)^2}{N} \frac{d\esp \left[Z_1^2 \norme{\bo{Z} }^{2} \right]}{d\esp \left[Z_1^2 \norme{\bo{Z} }^{2} \right] - 3 \rho_{\min}^2}  \left[\left( \frac{d\esp \left[Z_1^2 \norme{\bo{Z} }^{2} \right]}{3 \rho_{\min}^2}   \right)^{d_0}-1 \right] \, .   \nonumber        
\end{eqnarray}    

Finally, when $\esp \left[Z_1^2 \norme{\bo{Z} }^{2} \right] \geq 3 \rho_{\min}^2$ and 
knowing that  $\sum_{p=1}^{d_0} \binom{d}{p} \leq \sum_{p=1}^{d_0} d^p =\frac{d^{d_0}-1}{d-1} d \leq 2d^{d_0}$, we can write   
\begin{eqnarray}   
\var\left[\widehat{\M^c_{N, d_0}}(\bo{x}) \right]  & \leq &   \frac{\left(\digamma_{d_0}^{\max} \right)^2}{N}  \sum_{p=1}^{d_0}
\binom{d}{p}  \left( \frac{\esp \left[Z_1^2 \norme{\bo{Z}}^{2} \right]}{3 \rho_{\min}^2}   \right)^{p} \nonumber \\ 
 &\leq &   \frac{\left(\digamma_{d_0}^{\max} \right)^2}{N} \left(\frac{\esp \left[Z_1^2 \norme{\bo{Z}}^{2} \right]}{3 \rho_{\min}^2}   \right)^{d_0} \sum_{p=1}^{d_0} \binom{d}{p}
=   \frac{2 \left(\digamma_{d_0}^{\max} \right)^2}{N} \left(d\frac{\esp \left[Z_1^2 \norme{\bo{Z}}^{2} \right]}{3 \rho_{\min}^2}   \right)^{d_0}   \nonumber  \, .
\end{eqnarray}  
Also, note that $\sum_{p=1}^{d_0} \binom{d}{p} \leq 2^{d}$.     
           
\section{Proof of Corollary \ref{coro:mse1}} \label{app:coro:mse1}
Since $r^*=d_0-1$ and $|w|\leq d_0$, we have        
$
K_{1, \rho_{\min}, d_0-1} = \binom{d}{d_0} \left(\frac{1}{2  \rho_{\min}} \right)^{d_0}
$. 
The first result is obvious using the variance of the emulator provided in Theorem \ref{theo:rtestfc} and the bias from Corollary \ref{coro:biasfc1}.\\ 
The second result is due to the fact that  when  $G_j=F_j$, the terms in the Db-expansion of $\M$ are $L_2$-orthogonal.  
     
\section{Proof of Corollary \ref{coro:rate1}} \label{app:coro:rate1}
As $V_k \sim \mathcal{U}(-\xi, \xi)$, then $Z_k \sim \mathcal{U}(-\sqrt{3}, \sqrt{3})$ and  $\esp \left[Z_1^2 \norme{\bo{Z}}^{2} \right]= d+4/5$. Thus,  the first result holds.\\
For the second result, taking $\rho_{\min} = \sqrt{\frac{c_0d(d+0.8)}{3}}          
$ with $1 <c_0 \leq 2$ yields 
$$  
\var\left[\widehat{\M^c_{N, d_0}}(\bo{x}) \right]   \leq  \bar{\Sigma}_{d_0}=  \frac{\left(\digamma_{d_0}^{\max} \right)^2}{N} \left(\frac{1-(1/c_0)^{d_0}}{c_0-1}\right) \leq  \frac{\left(\digamma_{d_0}^{\max} \right)^2}{N (c_0-1)} \, .  
$$

\section{Proof of Corollary \ref{coro:mse2}} \label{app:coro:mse2} 
Using Equation (\ref{eq:intub}) and the fact that $r^*=d-1$, $d_0=d$, we can check that 
$$
\var\left[\widehat{\M^c_{N, d_0}}(\bo{x}) \right]   \leq   \frac{\left(\digamma_{d_0}^{\max} \right)^2}{N}  \left[\left( \frac{\esp \left[Z_1^2 \norme{\bo{Z} }^{2} \right]}{3 \rho_{\min}^2} +1  \right)^{d}-1 \right]
\leq  \bar{\Sigma}_{d}\, ,  
$$ 
with        
$$  
\bar{\Sigma}_{d} := \frac{d \left(\digamma_{d_0}^{\max} \right)^2 \esp \left[Z_1^2 \norme{\bo{Z} }^{2} \right]}{N\left(d\esp \left[Z_1^2 \norme{\bo{Z} }^{2} \right] - 3 \rho_{\min}^2\right)}  \left[\left( d\frac{\esp \left[Z_1^2 \norme{\bo{Z} }^{2} \right]}{3 \rho_{\min}^2}   \right)^{d}-1 \right]  \, .     
$$   
Thus, the results hold using Corollaries \ref{coro:biasfc2}, \ref{coro:mse1}.    

\end{appendices}                           
       


\end{document}